\documentclass[12pt]{article} 

\usepackage{amsmath,amsfonts,amssymb,amscd}

\def\dim{\mathop{\mathrm{dim}\;\!}\nolimits}
\def\codim{\mathop{\mathrm{codim}\;\!}\nolimits}

\def\rit{\mathbb{R}} 
\def\zit{\mathbb{Z}}   
\def\nit{\mathbb{N}} 
\def\cit{\mathbb{C}} 

\def\Cc{\mathcal{C}}
\def\Ac{\mathcal{A}}
\def\Lc{\mathcal{L}}
\def\Oc{\mathcal{O}}
\def\Ic{\mathcal{I}}
\def\S{\Sigma}
\def\tr{\mathrm{tr}}
\def\gg{\mathbf{g}}
\def\gh{\mathbf{h}}
\def\gnm{\mathbf{n_-}}
\def\gnp{\mathbf{n_+}}
\def\gb{\mathbf{b}}
\def\gn{\mathbf{n}}
\def\gb{\mathbf{b}}
\def\gbm{\mathbf{b_-}}
\def\gbp{\mathbf{b_+}}
\def\tnm{\mathbf{\widetilde{n_-}}}
\def\img{\mathrm{Im}}
\def\tf{\widetilde{f}}
\def\te{\widetilde{e}}
\def\thh{\widetilde{h}}
\def\tu{\widetilde{u}}
\def\tlam{\widetilde{\lambda}}
\def\talp{\widetilde{\alpha}}
\def\tLam{\widetilde{\Lambda}}
\def\tom{\widetilde{\omega}}
\def\tmu{\widetilde{\mu}}
\def\tgh{\widetilde{\gh}}
\def\bul{\bullet}

\def\la{\langle}
\def\ra{\rangle}

\def\tgg{\widetilde{\gg}}
\def\tM{\widetilde{M}}

\def\ftext#1{{\let\thefootnote\relax\footnotetext{\noindent #1}}}

\newcommand{\qed}{\hfill~~\mbox{$\square$}}

\newtheorem{theorem}{Theorem}[section]
\newtheorem{proposition}[theorem]{Proposition}

\newtheorem{lemma}[theorem]{Lemma}
\newtheorem{corollary}[theorem]{Corollary}

\begin{document}

\title{\Large {\bf Differential Equations Compatible with KZ Equations }}
\author{ \normalsize G. Felder$^*\,$, Y. Markov$^{\star}\,$, V. Tarasov$^{**}\,$
and A. Varchenko$^{\diamond}$}
\date{}
\maketitle
\begin{center}
{\it
$^*$ Departement Mathematik, ETH-Zentrum,\\
8092 Z\"urich, Switzerland\\ \hfill\\
$^{\star,\,\diamond}$ Department of Mathematics, University of North
Carolina,\\
Chapel Hill, NC 27599 -- 3250, USA\\ \hfill\\
$^{**}$ St. Petersburg Branch of Steklov Mathematical Institute\\
Fontanka 27, St. Petersburg, 191011, Russia}
\end{center}

\medskip
\centerline{January, 2000}
\medskip

\ftext{ \hspace{-0.6cm}
$^*$ {\sl E-mail\/{\rm:} felder@math.ethz.ch}\\
$^{\star}$ {\sl E-mail\/{\rm:}
markov@math.unc.edu}\\
$^{**}$ {\sl E-mail\/{\rm:} vt@pdmi.ras.ru}\\
$^{\diamond}$ {\sl E-mail\/{\rm:} av@math.unc.edu}\\
\hphantom{$*$} Supported in part by NSF grant DMS-9801582
}

\section{Introduction}

In the theory of the bispectral problem \cite{DG}, \cite{W}, one
considers a commutative algebra $A$ of differential operators
$L(z,\partial/\partial z)$ acting on functions of one complex variable $z$. Such an
algebra is called bispectral if there exists a non-trivial family
$u(z,\mu)$ of common eigenfunctions depending on a spectral parameter
$\mu$
\begin{equation}\label{firstequation}
Lu(z,\mu)=f_L(\mu)u(z,\mu), \qquad L\in A,
\end{equation}
which is also a family of common eigenfunctions of a commutative
algebra $B$ of differential operators $\Lambda(\mu,\partial/\partial\mu)$
with respect to $\mu$:
\begin{equation}\label{secondequation}
\Lambda u(z,\mu)=\theta_\Lambda(z)u(z,\mu), \qquad \Lambda\in B.
\end{equation}
J. Duistermaat and A. Gr\"unbaum \cite{DG} studied the case where $A$ is the algebra
of differential operators that commute with a Schr\"odinger operator
$\frac{d^2}{dz^2}-V(z)$ with meromorphic potential $V(z)$. They give a
complete classification of bispectral algebras arising in this way.
In particular they show
that $A$ is bispectral if  $V(z)$ is a rational KdV potential
(a rational function which stays rational under the flow
of the Korteweg--de Vries equation). G. Wilson \cite{W} classified bispectral
algebras of rank one, i.e., such that the greatest common divisor
of the orders of the differential operators in $A$ is one.
He showed that the maximal bispectral algebras of rank one
are in one to one correspondence with
conjugacy classes of
pairs $(Z,M)$ of square matrices so that $ZM-MZ+I$ has rank one.
The bispectrality then follows from the existence of the involution
$(Z,M)\mapsto (M^T,Z^T)$, which corresponds to exchanging $z$ and $\mu$.

The higher dimensional version of the bispectral problem, in which 
$A$ consists of partial differential operators in $z\in\cit^n$
is open. However O. Chalykh, M. Feigin and A. Veselov \cite{CV},\cite{CFV}
constructed examples of algebras in higher dimensions which have
the bispectral property (see Veselov's contribution to \cite{HK}). 
In these examples $A$  consists of differential operators commuting with an
$n$-particle Schr\"odinger operator with certain special rational potentials,
including the Calogero--Moser ones. These
potentials are in many respects the natural generalization of 
rational KdV potential associated to rank one algebras. In these
examples, the Baker--Akhiezer function $u(z,\mu)$ is symmetric
in the two arguments, thus $B=A$.

A good source of material on the bispectral problem is the volume
\cite{HK}.

In this paper we study 
a class of examples of commutative algebras of partial differential operators
acting on {\em vector-valued} functions with the bispectral property.
This means that in \eqref{firstequation}, \eqref{secondequation},
$u$ takes values in a vector space and $f_L(\mu)$, $\theta_\Lambda(z)$ 
are endomorphisms of the vector space.
In our class of examples, the algebra $A$ is generated by Knizhnik--Zamolodchikov
differential operators. They are commuting first order differential operators
associated to a complex simple Lie algebra
$\gg$ with a fixed non-degenerate invariant bilinear form and a non-zero complex
parameter $\kappa$. They act
on functions of $n$ complex variables $z_i$ taking values in the tensor
product of $n$ finite dimensional  $\gg$-modules. 
The ``dual'' variable $\mu$ is in a Cartan subalgebra of $\gg$. 
The first set of equations
\eqref{firstequation} is then  the set of generalized Knizhnik--Zamolodchikov
equations
\[
\left(\kappa\frac{\partial}{\partial z_i}-\sum_{j:j\neq i}\frac{\Omega^{(ij)}}{z_i-z_j}\right)u(z,\mu)
=\mu^{(i)}u(z,\mu), \qquad i=1,\dots,n.
\] 
Here $\Omega\in \gg\otimes\gg$ is dual to the invariant bilinear form and $\Omega^{(ij)}$ acts as
$\Omega$ on the $i$th and $j$th factors of the tensor product and as the identity on the other
factors. Similarly $\mu^{(i)}$ is $\mu$ acting on the $i$th factor. It is well-known
that these equations form a compatible system, i.e., they are the equations defining
horizontal sections for a flat connection. For $\mu=0$ they reduce to the classical
Knizhnik--Zamolodchikov equations.
The algebra $B$ is generated by rank($\gg$) first order partial differential operators in $\mu$
with rational coefficients. We call the corresponding equations \eqref{secondequation}
{\em dynamical differential equations}, and show that they form, together with the 
generalized Knizhnik--Zamolodchikov equations, a compatible system.
We also give simultaneous solutions of both systems of equations
in terms of hypergeometric integrals for a more  general class of Lie algebras, 
which includes in particular all Kac--Moody Lie algebras.

In the case of $\gg=sl_2$, the algebra $B$ is generated by one ordinary differential
operator. In this case, the corresponding equations where first written and solved by
H. Babujian and A. Kitaev \cite{BK}, who also related the equations to the Maxwell--Bloch
system.

Our paper is organized as follows. In Sect.\ \ref{dualKZ} we introduce
the systems of Knizhnik--Zamolodchikov and dynamical differential equations for
arbitrary simple Lie algebra and prove their compatibility.  We then give
formulae for hypergeometric solutions in Sect.\ \ref{HG-sol} and give as an application a
determinant formula. The fact that the hypergeometric integrals provide
solutions is a consequence of a general theorem valid for a class of Lie
algebras with generic Cartan matrix, introduced in \cite{SV1}. We introduce
in Sect.\ \ref{free}
the Knizhnik--Zamolodchikov and dynamical differential equations in this
more general context and explain in the next Section the results on
complexes of hypergeometric differential forms from \cite{SV1}.
In Sect.\ \ref{DynDE} we prove that the hypergeometric integrals
for generic Lie algebras satisfies the dynamical differential equations. 
Finally in Sect.\ \ref{main-theorems}
we prove that our hypergeometric integrals are solutions of both
systems of equations for any Kac--Moody Lie algebra. We also find
a determinant formula, which implies a completeness result for solutions
in the case of generic parameters.

\section{Dynamical differential equations}\label{dualKZ}
\subsection{}
  Let $\gg$ be a simple complex Lie algebra with an invariant bilinear
form $(\, , \, )$ and a root space decomposition $\gg =
\gh\oplus(\oplus_{\alpha\in\Delta}\cit e_{\alpha})$. The
root vectors $e_{\alpha}$ are normalized so that
$(e_{\alpha},e_{-\alpha})=1$. Then the quadratic Casimir element of
$\gg\otimes \gg$ has the form $\Omega = \sum_s h_s\otimes h_s +
\sum_{\alpha \in \Delta} e_{\alpha}\otimes e_{-\alpha}$, for any
orthonormal basis
$(h_s)$ of the Cartan subalgebra $\gh$. We also fix a system of
simple roots $\alpha_1,\ldots ,\alpha_r$.

Consider the Knizhnik--Zamolodchikov (KZ) equations with an additional parameter $\mu \in \gh$,
for a function $u$ on $n$ variables taking values in a tensor product
$V=V_1 \otimes\cdots\otimes V_n$ of highest weight modules of $\gg$ with
corresponding highest weights $\Lambda_1,\ldots,\Lambda_n$,
\begin{equation}\label{KZ}
 \kappa\frac{\partial u}{\partial z_i}=\mu^{(i)}u + \sum_{i\ne j}
 \frac{\Omega^{(ij)}}{z_i - z_j}u, \qquad\qquad i=1,\ldots,n,
\end{equation}
where $\kappa$ is a complex parameter.
 We are interested in a differential equation for $u$ with respect to
$\mu$ which are compatible with KZ equations. If $\mu' \in \gh$,
denote by $\partial_{\mu'}$ the partial derivative with respect to
$\mu$ in the
direction of $\mu'$
\begin{theorem}\label{compatibility} The equations
\begin{equation}\label{KZd}
 \kappa\partial_{\mu'}u=
\sum_{i= 1}^{n}z_i(\mu')^{(i)}u + 
\sum_{\alpha>0
}\frac{\la\alpha,\mu'\ra}{\la\alpha,\mu\ra}e_{-\alpha}e_{\alpha}u,
 \qquad\qquad \mu' \in \gh
\end{equation}
form together with the KZ equations (\ref{KZ}), a compatible 
system of equations for a function $u(z,\mu)$ taking values in 
$V=V_1 \otimes\cdots\otimes V_n$.
\end{theorem}
 The equations (\ref{KZd}) will be called {\it dynamical differential equations}.

\medskip
\noindent{\it{Example}}. Let $\gg=sl_N= gl_N/\cit$. View $sl_N$-modules as $gl_N$-modules
by letting the center of $gl_N$ act trivially. Denote by $E_{a,b}\in gl_N$ the
matrix whose entries are zero except for a one at the intersection of the $a$th row with
the $b$th column. The fundamental coweights $\varpi_a=(1-a/N)\sum_{b\leq a}E_{b,b}-(a/N)\sum_{b>a}E_{b,b}$,
$a=1,\dots,N-1$ form a basis of the standard Cartan subalgebra of $sl_N$. Write 
$\mu=\sum_{a=1}^{N-1}\mu_a\varpi_a$.
 Then our equations may be written as
\begin{eqnarray*}
 \kappa\frac{\partial u}{\partial z_i}  & = &
\sum_{a=1}^{N-1}\mu_{a}\varpi_a^{(i)}u + \sum_{j:j\ne i}
\sum_{a,b}\frac{E_{a,b}^{(i)}E_{b,a}^{(j)}}{z_i - z_j}\,u.\\
 \kappa\frac{\partial u}{\partial \mu_a}  & = &
\sum_{i=1}^n z_i\varpi_a^{(i)}u + \sum_{b,c:b\leq a<c}
\sum_{i,j}\frac{E_{b,c}^{(i)}E_{c,b}^{(j)}}{\mu_b+\mu_{b+1}\cdots+\mu_{c-1}}\,u.
\end{eqnarray*}

\subsection{Proof of Theorem~\ref{compatibility}}
It is rather easy to verify that most terms of the compatibility
equations vanish. The only non-trivial thing to check is that the
operators $\sum_{\alpha >0}
\frac{\la\alpha,\lambda\ra}{\la\alpha,\mu\ra}e_{-\alpha}e_{\alpha}$ commute
for different values of $\lambda$. The operators obtained by
extending the sum to all roots differ from the sum over
positive roots by an element of the Cartan subalgebra.
Since the operators commute with the Cartan subalgebra, it
is sufficient to prove the following Proposition.

\begin{proposition}\label{t-commute}
 Let for $\lambda, \mu \in \gh$, $T(\lambda,\mu)=
\sum_{\alpha\in\Delta}
\frac{\la\alpha,\lambda\ra}{\la\alpha,\mu\ra}e_{-\alpha}e_{\alpha}$. Then for
all $\lambda,\mu\in\gh$, 
$$T(\lambda,\mu)T(\nu,\mu)=T(\nu,\mu)T(\lambda,\mu).$$
\end{proposition}
  The proof is based on the following fact.
\begin{lemma} 
 Let $\alpha, \beta \in \Delta$ with $\alpha \ne \pm\beta$, and let
$S=S(\alpha,\beta)$ be the set of integers $j$ such that
$\beta+j\alpha\in\Delta$. Then $\sum_{j\in S}
[e_{\alpha},e_{\beta +j\alpha}e_{-\beta -j\alpha}]=0$ and 
$\sum_{j\in S}
[e_{-\alpha},e_{\beta +j\alpha}e_{-\beta -j\alpha}]=0$
\end{lemma}
{\it{Proof.}} For roots $\gamma,\delta$ such that $\gamma+ \delta$ is
a root, let $N_{\gamma,\delta}=([e_{\gamma},e_{\delta}],
e_{-\gamma-\delta})$, so that
$[e_{\gamma},e_{\delta}]=N_{\gamma,\delta}e_{\gamma+ \delta}$. By
considering the adjoint action on $\gg$ of the $sl_2$ sub-algebra
generated by $e_{\pm\alpha}$, we see that for $\beta \ne
\pm\alpha$, $S$ is a finite sequence of subsequent integers. We may
thus assume that $S=\{ 0, \ldots, k\}$ by replacing $\beta$ by
$\beta-j\alpha$ for some $j$ if necessary. 

We then have
\begin{eqnarray*}
 \lefteqn{\sum_{j=0}^k
[e_{\alpha},e_{\beta +j\alpha}e_{-\beta-j\alpha}]} &&\\
 & = & \sum_{j=0}^{k-1}N_{\alpha,\beta +j\alpha}e_{\beta
+(j+1)\alpha}e_{-\beta-j\alpha}+
 \sum_{j=1}^{k}N_{\alpha,-\beta -j\alpha}e_{\beta
+j\alpha}e_{-\beta-(j-1)\alpha}\\
 & = &  \sum_{j=0}^{k-1}(N_{\alpha,\beta +j\alpha}+N_{\alpha,-\beta
-(j+1)\alpha})e_{\beta+(j+1)\alpha}e_{-\beta-j\alpha}
\end{eqnarray*}
By the invariance of the bilinear form,
\begin{eqnarray*}
N_{\alpha,\beta +j\alpha} & = & 
  ([e_{\alpha},e_{\beta+j\alpha}], e_{-\beta-(j+1)\alpha})\\
 & = & -(e_{\beta+j\alpha},[e_{\alpha},e_{-\beta-(j+1)\alpha}])\\
 & = & -N_{\alpha,-\beta-(j+1)\alpha}. 
\end{eqnarray*}
Therefore $\sum_{j\in S}
[e_{\alpha},e_{\beta +j\alpha}e_{-\beta -j\alpha}]$ vanishes. The
other statement is proved by replacing $\alpha$ by $-\alpha$ and
noticing that $S(-\alpha,\beta)=S(\alpha,\beta)$. \qed

{\it{Proof of Proposition~\ref{t-commute}:}} Consider
\begin{equation}\label{for-t-comm}
T(\lambda,\mu)T(\nu,\mu)=\sum_{\alpha,\beta}\frac{\la\alpha,\lambda\ra
\la\beta,\nu\ra}{\la\alpha,\mu\ra\la\alpha,\mu\ra}
[e_{\alpha}e_{-\alpha},e_{\beta}e_{-\beta}]
\end{equation}
Let us show that this expression is a regular function of
$\mu\in\gh$. Since for nontrivial $\lambda, \mu$ it converges to zero
at infinity, it then vanishes identically.

We compute the residue of (\ref{for-t-comm}) at $\la\alpha,\mu\ra=0$:
\begin{equation}\label{sum-gamma}
\sum_{\gamma: \gamma\ne\pm\alpha}\frac{\la\alpha,\lambda\ra 
\la\gamma,\nu\ra- \la\gamma,\lambda\ra\la\alpha,\nu\ra}{\la\gamma,\mu\ra}
[e_{\alpha}e_{-\alpha},e_{\gamma}e_{-\gamma}],
\end{equation}
a function on the hyperplane  $\la\alpha,\mu\ra=0$. The sum over $\gamma$
of the form $\beta +j\alpha$, $j\in S(\alpha,\beta)$ gives for
 $\la\alpha,\mu\ra=0$,
\begin{equation}
\sum_{j\in S}\frac{\la\alpha,\lambda\ra\la\beta,\nu\ra-
 \la\beta,\lambda\ra\la\alpha,\nu\ra}{\la\beta,\mu\ra} 
[e_{\alpha}e_{-\alpha},e_{\beta+j\alpha}e_{-\beta-j\alpha}]=0,
\end{equation}
by the previous Lemma. Since the sum over $\gamma$ in
(\ref{sum-gamma}) can be written as a sum of such terms, it vanishes.

\section{Hypergeometric solutions}\label{HG-sol}

Let $\gg$ be a simple complex Lie algebra. Choose a set $f_1,\ldots,f_r$, 
$e_1,\ldots,e_r$ of Chevalley generators of the Lie algebra $\gg$ associated 
with simple roots $\alpha_1,\ldots,\alpha_r$. Let $\lambda=(m_1,\ldots,m_r)\in\nit^r$.
Let $Q_+=\sum \nit\alpha_i$ be the {\it positive root lattice} for $\gg$. Define a map 
$\alpha:\nit^r\to Q_+$ by $\alpha(\lambda)=\sum m_i\alpha_i$. Let $V$ be a tensor product of
highest weight modules $V_j$ of $\gg$ with respective highest weights $\Lambda_j$,
where $j=1,\ldots,n$. Set $\Lambda=\sum_{j=1}^n\Lambda_j$. Denote $V_{\lambda}$ the weight space
of $V$ with weight $\Lambda - \alpha(\lambda)$.
 The hypergeometric solutions of the KZ equations in $V_{\lambda}$, see \cite{SV1}, 
have the form  
$u(z)=\int_{\gamma(z)}\Phi(z,t)^{1/k}\omega(z,t)$. We will describe the explicit construction.
The number of integration 
variables $(t_k)_{k=1}^m$ is $m=\sum m_i$.
Let $c$ be the unique non-decreasing function from 
$\{1,\ldots,m\}$ to $\{1,\ldots,r\}$ $(i=1,\ldots,r)$, such that 
$\# c^{-1}(\{i\})=m_i$. Define
$$ \Phi(z,t)=\prod_{i<j}(z_i-z_j)^{(\Lambda_i,\Lambda_j)}
\prod_{k,j}(t_k-z_j)^{-(\alpha_{c(k)},\Lambda_j)}
\prod_{k<l}(t_k-t_l)^{(\alpha_{c(k)},\alpha_{c(l)})}.
$$
 The $m$--form $\omega(z,t)$ is a closed logarithmic 
differential form on $\cit^n\times\cit^m$ with values in $V_{\lambda}$.
It has the following combinatorial description.  Let $P(\lambda,n)$ be the set of sequences 
$I=(i_1^1,\ldots,i^1_{s_1};\ldots;i^n_1,\ldots,i^n_{s_n})$ of integers in 
$\{1,\ldots,r\}$ with $s_j\geq 0$, $j=1,\ldots,n$ and such that, for all
$1\leq j\leq r$, $j$ appears precisely $|c^{-1}(j)|$ times in $I$. For 
$I\in P(\lambda,n)$, and a permutation $\sigma\in \Sigma_m$, set $\sigma_1(l)=\sigma(l)$ 
and $\sigma_j(l)=\sigma(s_1+\cdots +s_{j-1}+l)$, $j=2,\ldots,n$, 
$1\leq l\leq s_j$. Define $\Sigma(I)=\{\sigma\in \Sigma_m\, |\, c(\sigma_j(l))=i_j^l$ 
for all $j$ and $l\}.$ 

Fix a highest-weight vector $v_j$ for each representation $V_j$, 
$j=1,\ldots,n$. To every $I\in P(\lambda,n)$ we associate a vector 
$f_Iv=f_{i_1^1}\cdots f_{i_{s_1}^1}v_1\otimes\cdots\otimes f_{i_1^n}\cdots 
f_{i_{s_n}^n}v_n$ in $V_{\lambda}$, and meromorphic differential $m$--forms 
$\omega_{I,\sigma}=\omega_{\sigma_1(1),\ldots,\sigma_1(s_1)}(z_1)\wedge\cdots
\wedge\omega_{\sigma_n(1),\ldots,\sigma_n(s_n)}(z_n)$,
 labeled by $\sigma\in \Sigma(I)$,
where $\omega_{i_1,\ldots, i_s}(z)=d\log(t_{i_1}-t_{i_2})\wedge\cdots\wedge
d\log(t_{i_{s-1}}-t_{i_s})\wedge d\log(t_{i_s}-z)$ is a meromorphic one form on
$\cit\times\cit^{s}$. Finally 
$$\omega(z,t)=\sum_{I\in P(\lambda,n)}\sum_{\sigma\in \S(I)}(-1)^{|\sigma|}
\omega_{I,\sigma}f_Iv.$$

 It obeys the equation
$\Phi^{-1}d\Phi\wedge\omega = K\wedge\omega$, where 
$K=\sum_{i<j}\Omega^{(ij)}\frac{dz_i-dz_j}{z_i-z_j}$. As a
consequence, for each horizontal family of twisted cycles $\gamma(z)$
in $\{ z\}\times\cit^m$, $u$ obeys the KZ equations (\ref{KZ}) with
$\mu=0$ (See~\cite{SV1}). This construction can be modified to give
solutions for general $\mu$:
\begin{theorem}\label{main1}
The integrals
\begin{equation}
\int_{\gamma(z)} \Phi_{\mu}^{\frac{1}{k}}\omega,\qquad\qquad
\Phi_{\mu}=\exp (-\sum_{i=1}^m\la\alpha_{c(i)},\mu\ra t_i +
\sum_{j=1}^n\la\Lambda_j,\mu\ra z_j)\Phi ,
\end{equation}
are solutions of the KZ equations (\ref{KZ}).
\end{theorem}
{\it{Proof.}} The proof follows from the identity
$$\Phi_{\mu}^{-1}d\Phi_{\mu}\wedge\omega = 
\sum_{i=1}^n\mu^{(i)}dz_i\wedge\omega + K\wedge\omega.$$ 
To prove this identity  notice that $\omega$ is a sum of several terms 
$\omega_s$, which can be grouped according to how the integration variables 
$t_1,\ldots, t_m$ are distributed among the points $z_1,\ldots, z_n$. If the 
variable $t_k$ is associated to the point $z_i$ in one term $\omega_s$, then
$(dt_k-dz_i)\wedge\omega_s=0$. Moreover, if the term $\omega_s$ obeys 
$h^{(i)}\omega_s=(\la\Lambda_i,h\ra -\sum m'_j(i)\la\alpha_j,h\ra )\omega_s$, $h\in\gh$, 
where $\Lambda_i$ is the highest weight of the {\it{i}}-th representation,
then the number of variables $t_k$ associated to $z_i$ such that $c(k)=j$ is
$m'_j(i)$, $(j=1,\ldots,r)$. Thus we have
$$\Phi_{\mu}^{-1}d\Phi_{\mu}\wedge\omega = 
\sum_{i=1}^n (\la\Lambda_i,\mu\ra -\sum_{j} m'_j(i)\la\alpha_j,\mu\ra )
dz_i\wedge\omega + K\wedge\omega.$$ The proof is complete.
\qed
\begin{theorem}\label{main2}
The hypergeometric integrals of Theorem~\ref{main1} obey the Dynamical differential
equations (\ref{KZd}).
\end{theorem}
The proof of the Theorem is given in Section~\ref{main-theorems} .
\qed

An application of the above two theorems is a determinant formula.
\begin{corollary} Fix a basis $v_1,\ldots, v_d$ of a weight space $V_{\lambda}$. Suppose that
$u_i(\mu,z)=\sum_{j=1}^d u_{i,j}v_j$, $i=1,\ldots,d$ is a basis of the space of solutions in a 
neighbourhood of a generic point $(\mu,z)\in\gh\times\cit^n$. Let $\delta_{\alpha}=
\tr_{V_{\lambda}}(e_{-\alpha}e_{\alpha})$, ($\alpha\in\Delta,\, \alpha>0$), $\epsilon_{ij}=
\tr_{V_{\lambda}}(\Omega_{ij})$. Then there is a constant $C=C(V_1,\ldots,V_n,\lambda,\kappa)
\ne 0$ such that
$$ \det (u_{ij})= C\exp(\sum_{i=1}^n \frac{z_i}{\kappa}\tr_{V_{\lambda}}(\mu^{(i)}))
\prod_{\alpha>0}\la\alpha,\mu\ra^{\frac{\delta_{\alpha}}{\kappa}}
\prod_{i<j}(z_i-z_j)^{\frac{\epsilon_{ij}}{\kappa}}.$$
\end{corollary}

\section{Free Lie algebras.  Dynamical differential equations}\label{free}
\subsection{The definition of KZ and dynamical differential
equations.}\label{defKZ}
 Following \cite{SV1} let us fix the following data:

 1. A finite dimensional complex vector space $\gh$;

 2. A non-degenerate symmetric bilinear form $(\, , \, )$ on $\gh$;

 3. Linearly independent covectors (``simple roots'') 
$\alpha_1,\ldots,\alpha_r\in\gh^{\ast}$.\\
 We denote by $b:\gh\rightarrow\gh^{\ast}$ the isomorphism induced by
$(\, , \, )$, and we transfer the form $(\, , \, )$ to $\gh^{\ast}$ via $b$.
Set $b_{ij}=(\alpha_i,\alpha_j)$, $h_i=b^{-1}(\alpha_i)\in\gh$. Denote by 
$\gg$ the Lie algebra generated by $e_i$, $f_i$ for $i=1,\ldots,r$ and $\gh$
subject to the relations:
$$ [h,e_i]=\la\alpha_i,h\ra e_i, \quad [h,f_i]=-\la\alpha_i,h\ra f_i,\quad
[e_i,f_j]=\delta_{ij}h_i,\quad [h,h']=0,$$
for all $i,j=1,\ldots,r$ and $h,h'\in \gh$. Thus we have constructed a 
Kac-Moody Lie algebra without Serre's relations. We denote by $\gnm$
(resp. by $\gnp$) the subalgebra of $\gg$ generated by $f_i$ (resp. 
$e_i$) for $i=1,\ldots,r$. We have $\gg=\gnm\oplus\gh\oplus\gnp$. Set 
$\gb_{\pm}=\gn_{\pm}\oplus\gh$. These are subalgebras of $\gg$.

Let $\Lambda\in\gh^{\ast}$. Denote by $M(\Lambda)$ the Verma module over $\gg$
generated by a vector $v$ subject to the relations $\gnp v=0$ and 
$hv=\la\Lambda,h\ra v$ for all $h\in\gh$. Let us fix weights $\Lambda_1,
\ldots,\Lambda_n\in \gh^{\ast}$, and let 
$M=M(\Lambda_1)\otimes \cdots\otimes M(\Lambda_n)$.

 For $\lambda =(m_1,\ldots,m_r)\in \nit^r$, set
\begin{eqnarray*} 
(\gn_{\pm})_{\lambda} & = &\left\{x\in\gn_{\pm}\,|
\, [h,x]=\la\pm\sum m_i\alpha_i,h\ra x, \quad \mbox{ for all } h\in\gh\right\},\\
M(\Lambda)_{\lambda} & = & \left\{x\in M(\Lambda)\,|\, hx=
\la\Lambda-\sum m_i\alpha_i,h\ra x \quad \mbox{ for all } h\in\gh\right\},\\
M_{\lambda}  & = & \left\{x\in M\,|\, hx=
\la\sum\Lambda_j-\sum m_i\alpha_i,h\ra x \quad \mbox{ for all } h\in\gh\right\}.
\end{eqnarray*}
We have $\gn_{\pm}=\oplus_{\lambda}(\gn_{\pm})_{\lambda}$, 
 $\gb_{\pm}=\gh\oplus (\oplus_{\lambda}(\gn_{\pm})_{\lambda})$, 
$M(\Lambda)=\oplus_{\lambda}M(\Lambda)_{\lambda}$,
$M=\oplus_{\lambda}M_{\lambda}$.

 Let $\tau : \gg \rightarrow \gg$ be the Lie algebra automorphism such that
$\tau(e_i)=-f_i$, $\tau(f_i)=-e_i$, $\tau(h)=-h$, for $h\in\gh$. 
Set $\gn_{\pm}^*=\oplus_{\lambda}(\gn_{\pm})_{\lambda}^*$. Set 
$M(\Lambda)^{\ast}=\oplus_{\lambda}M(\Lambda)_{\lambda}^{\ast}$. Define a
 structure of  a $\gg$-module on $M(\Lambda)^{\ast}$ by the rule
\begin{equation}\label{*action}
\la g\phi,x\ra=\la\phi,-\tau(g)x\ra \mbox{ for } \phi \in M(\Lambda)^{\ast},\quad g\in\gg, 
\quad x\in M(\Lambda).
\end{equation}
There is a unique bilinear form $K(\, ,\,)$ on $\gg$ such that: 
$K$ coincides with $(\, ,\,)$ on $\gh$; $K$ is zero on $\gnp$ and $\gnm$;
$\gh$ and $\gnm\oplus\gnp$ are orthogonal; 
$K(f_i,e_j)=K(e_j,f_i)=\delta_{ij}$ for $i,j=1,\ldots,r$;
$K$ is $\gg$-invariant, that is $K([x,y],z)=K(x,[y,z])$ for all $x,y,z\in\gg$.

 A bilinear form $S$ on $\gg$ is defined by the rule $S(x,y)=-K(\tau(x),y)$.
The form $S$ is symmetric, $\tau$-invariant, and $S([x,y],z)=S(x,[\tau(y),z])$.
The subspaces $\gnp$, $\gh$, $\gnm$ are pairwise orthogonal with respect to 
$S$.

 For a Verma module $M(\Lambda)$ with highest weight $\Lambda$ and generating
vector $v$ there is  a unique bilinear form $S$ on $M(\Lambda)$ such that:
$S(v,v)=1$; \space $S(e_ix,y)=S(x,f_iy)$; \space
$S(f_ix,y)=S(x,e_iy)$, for all $x,y\in M(\Lambda)$ and $i=1,\ldots,r$.
$S$ is symmetric. The subspaces $M(\Lambda)_{\lambda}$ are pairwise 
orthogonal with respect to $S$. The form $S$ induces a homomorphism of
$\gg$-modules $S:M(\Lambda)\rightarrow M(\Lambda)^{\ast}$. The module 
$M(\Lambda)/ \ker S$ is the irreducible $\gg$-module with highest weight
$\Lambda$. More generally, on the space $\wedge^p\gnm\otimes M(\Lambda_1)
\otimes\cdots\otimes M(\Lambda_n)$ for $\Lambda_j\in\gh^{\ast}$, a bilinear
form S is defined by the rule
$$ S(g_1\wedge\cdots\wedge g_p\otimes x_1\otimes\cdots\otimes x_n,
g'_1\wedge\cdots\wedge g'_p\otimes x'_1\otimes\cdots\otimes x'_n)=
\det S(g_i,g'_j)\cdot\prod_{i=1}^n S(x_i,x'_i). $$
This form induces a map
\begin{equation}\label{contrav}M
 S:\wedge^p\gnm\otimes M(\Lambda_1)\otimes\cdots\otimes M(\Lambda_n)
\rightarrow (\wedge^p\gnm\otimes M(\Lambda_1)\otimes\cdots\otimes 
M(\Lambda_n))^{\ast}
\end{equation}
$S$ is called the contravariant form. The linear map (\ref{contrav})
depends analytically on $(b_{ij})$ and 
$\Lambda_1,\ldots,\Lambda_n\in \gh^{\ast}$ for a fixed $r$. It is 
non-degenerate
for general values of the parameters, see Theorem~3.7 in \cite{SV1}.

 {\bf{A Lie bialgebra structure on $\gb_{\pm}$ \cite{SV1}.}}
A Lie bialgebra is a vector space $\gg$ with a Lie algebra structure and a Lie coalgebra
structure, such that the cocommutator map $\nu :\gg\rightarrow \gg\wedge\gg$ 
is a  one-cocycle: $x\nu(y) - y\nu(x)=\nu([x,y])$, for all $x,y\in\gg$. Here the 
action of $\gg$ on $\gg\wedge\gg$ is the adjoint one: 
$a(b\wedge c)=[a,b]\wedge c + b\wedge [a,c]$. The dual map to the cocommutator map,
$\nu^*:(\gg\wedge\gg)^*\rightarrow \gg^*$,
defines a Lie algebra structure on $\gg^{\ast}$.

 Let $\gg$ be a Lie bialgebra. The double of $\gg$ is the Lie algebra equal to
$\gg\oplus\gg^{\ast}$ as a vector space with the bracket on $\gg$ and $\gg^{\ast}$ defined by
the Lie algebra structure on $\gg$ and $\gg^{\ast}$, and for $x\in\gg$ and 
$l\in\gg^{\ast}$, $[l,x]=\bar{l}+\bar{x}$, where $\bar{x}\in\gg$ and 
$\bar{l}\in\gg^{\ast}$ are defined by the rules $\bar{l}(y)=l([x,y])$, 
$m(\bar{x})=[m,l](x)$, for all $y\in\gg$, $m\in\gg^{\ast}$. The double is denoted by
$D(\gg)$.

 Let $\gg$ be 
the Kac-Moody algebra we defined at the beginning of the section. $\gg$ is
a Lie bialgebra with respect to the following cobracket. There exists a
 unique map $\nu : \gg\rightarrow\gg\wedge\gg$ such that 
$x\nu(y)-y\nu(x)=\nu([x,y])$, and $\nu(h)=0$, and 
$\nu(f_i)=\displaystyle\frac12f_i\wedge h_i$, and
$\nu(e_i)=\displaystyle\frac12e_i\wedge h_i$. In the previous four equalities,
$h\in\gh$, $i=1,\ldots,r$, the action of $\gg$ on $\gg\wedge\gg$ is the adjoint one,
see \cite{D1} Example 3.2 and \cite{SV1}.
$\gbm$ and $\gbp$ are subbialgebras . The map $\nu$ has the property 
$\tau\nu +\nu\tau=0$. Thus, if $\rho :\gbm^{\ast}\rightarrow\mathrm{End}(V)$
is a representation of the Lie algebra 
$(\nu |_{\gbm})^{\ast}:\Lambda^2\gbm^{\ast}\rightarrow \gbm^{\ast}$, then
 $-\rho\circ\tau :\gbp^{\ast}\rightarrow\mathrm{End}(V)$
is a representation of the Lie algebra 
$(\nu |_{\gbp})^{\ast}:\Lambda^2\gbp^{\ast}\rightarrow \gbp^{\ast}$.

  Note that the coalgebra map $\nu$ defines a Lie algebra structure on
$\gb_{\pm}^{\ast}$.

 {\bf{Comultiplication.}} Let $M=M(\Lambda_1)\otimes\cdots\otimes 
M(\Lambda_n)$.
Let $v=v_1\otimes\cdots\otimes v_n\in M$ be the product of the generating
vectors. Set $\Lambda =\sum\Lambda_j$. Let $\gbm$ act on $\gbm\otimes M$ by 
the rule $a(b\otimes m)=
[a,b]\otimes m + b\otimes am$. For $1\leq i\leq n$, and $a,b\in\gg$, 
$m=x_1\otimes\cdots\otimes x_n\in M$, set $a^{(i)}m=x_1\otimes\cdots\otimes 
x_{i-1}\otimes ax_i\otimes x_{i+1}\otimes\cdots\otimes x_n$ and 
$a^{(i)}(b\otimes m)=[a,b]\otimes m + b\otimes a^{(i)}m$.
 
  There is a unique linear map $\nu_M: M\rightarrow \gbm\otimes M$ such that 
\begin{align}
&\nu_M(h\cdot x)=h\cdot\nu_M(x) \mbox{ for any } h\in\gh \mbox{ and } x\in M;\notag\\ 
&\nu_M(x)=\frac12(b^{-1}(\Lambda-\alpha(\lambda)))\otimes x  +
\nu_{M-}(x) \mbox{ for } x\in M_{\lambda}.\notag
\end{align}
Recall that $b^{-1}(\alpha(\lambda))=\sum m_ib^{-1}(\alpha_i)=\sum m_ih_i$,
and $b^{-1}:\gh^{\ast}\rightarrow\gh$ is defined at the beginning of Sect.~\ref{defKZ}.
The map $\nu_{M-}: M\to\gnm\otimes M$ is defined via an inductive definition.
$\nu_{M-}(v)=0$, $\nu_{M-}(x)=\sum_{k=1}^n\nu_{M-}^{(k)}(x)$ where 
$$\nu^{(k)}_{M-}(f_i^{(j)}x)
=f_i^{(j)}\nu^{(k)}_{M-}(x) \mbox{ for  } k\ne j, 
\mbox{ and }
\nu^{(k)}_{M-}(f_i^{(k)}x)=f_i\otimes h_i^{(k)}x +
f_i^{(k)}\nu^{(k)}_{M-}(x).$$ 
In all formulae we have $1\leq i\leq r$, and $1\leq j,k\leq n$.

 {\bf{Remark.}} The corresponding definition of $\nu_{M-}(x)$ in
\cite{SV1}
should be corrected as above. Note that the two definitions coincide if we have one 
tensor factor, i.e. $n=1$.
 We have the following Lemma (cf. Lemma~6.15.2 in \cite{SV1}).

\begin{lemma}\label{CDlemma} For any $x,y\in M$, $a\in \gbm$,
\begin{equation}
 S(\nu_M(x),a\otimes y)=
 \begin{cases}
  \frac12S(x,ay) & \text{ if } a\in\gh; \\
  S(x,ay)        & \text{ if } a\in\gnm.
 \end{cases}
\end{equation}
Here $S$ is defined on $\gbm\otimes M$ by the rule
$S(a\otimes x, b\otimes y)=S(a,b)S(x,y)$, cf. (\ref{contrav}).
\end{lemma}
The proof of Lemma~\ref{CDlemma} is given in Section~\ref{proofCDlem} \qed.\\
Note that the above Lemma renders the following diagram commutative.
\[ \begin{CD}
    \gnm\otimes M @>{standard}>> M\\
    @VV{S}V   @VV{S}V\\
    \gnm^{\ast}\otimes M^{\ast} @>{\nu_{M-}^{\ast}}>> M^{\ast}
\end{CD} \]
 $S$ is an isomorphism for general values of parameters $(b_{ij})$, 
$(\Lambda_k)_{k=1}^n$, see \cite{SV1} sec. (3.7) and (6.6). Hence,
$\nu_M^{\ast}:\gbm^{\ast}\otimes M^{\ast}\rightarrow M^{\ast}$ is a
$\gbm^{\ast}$-module structure with respect to the Lie algebra structure
$\nu^{\ast}:\Lambda^2\gbm^{\ast}\rightarrow\gbm^{\ast}$ for any values of the above 
parameters.
\begin{corollary}\label{CDcor}
 If $y\in \ker S:\gnm\rightarrow\gnm$ then $y M\subset \ker S:M
\rightarrow M$. \qed\\
\end{corollary}

  {\bf{Actions of the doubles of $\gb_{\pm}$ on $M$ and $M^{\ast}$ 
respectively.}}

  Consider the standard action of $\gbm$ on $M^{\ast}$, i.e. $\forall a\in\gbm$ $\forall
\phi\in M^*$, $\la a\cdot\phi,\,.\,\ra= \la \phi,-a\,.\,\ra$ 
where the action on the right hand side 
is the standard action of $\gg$ on $M$. This map together with
$\nu_{M-}^{\ast}$ defines an action of $\gbm\oplus\gbm^{\ast}$ on $M^{\ast}$.
Lemma~6.17.1 \cite{SV1} asserts that $M^{\ast}$ is a $D(\gbm)$-module
under this action, where  $D(\gbm)$ is the double of $\gbm$.


  For any $a\in \gbm^{\ast}$, the action of $\nu_M^{\ast}$ defines a map
$\nu_M^{\ast}(a,\, .):M^{\ast}\rightarrow M^{\ast}$. Set 
$\rho(a)=-(\nu_M^{\ast}(a,\, .))^{\ast}:M\rightarrow M$. The map
$\rho:\gbm^{\ast}\rightarrow\mathrm{End}(M)$ gives an action of 
$\gbm^{\ast}$ on $M$. An action of $\gbp^{\ast}$ on $M$ is defined by
$\omega =-\rho\circ\tau:\gbp^{\ast}\rightarrow\mathrm{End}(M)$. The rule
$a\otimes x\rightarrow \tau(a)x$, for $a\in\gbp, x\in M$, defines an action 
of $\gbp$ on $M$. This action and the map $\omega$ define an action of 
$\gbp\otimes\gbp^{\ast}$ on $M$. Lemma~11.3.28 \cite{V1} implies that $M$ is 
a $D(\gbp)$ module.

  Note that for any Kac-Moody Lie algebra $\gg$ without Serre's relations, 
and a $\gg$-module $V$ we can define a $\gg$-module structure on $V^{\ast}$
by the rule $\la g\phi,\,\cdot\,\ra=\la\phi,-g\,\cdot\,\ra$ for all $g\in \gg$ and 
$\phi\in M^{\ast}$. With respect to this module structure, $M$ becomes a
$D(\gbm)$-module and $M^{\ast}$ becomes a $D(\gbp)$-module.

{\bf{KZ equations and dynamical differential equations in $M$ and
$M^{\ast}$.}} For 
a vector space $V$ denote by $\Omega(V)\in V\otimes V^{\ast}$ the canonical 
element. For $\lambda\in\nit^r$, set $\Omega_{\lambda,\pm}^-:=
\Omega((\gn_{\pm})_{\lambda}) \in (\gn_{\pm})_{\lambda}\otimes
(\gn_{\pm})_{\lambda}^{\ast}$, 
 set $\Omega_{\lambda,\pm}^+:=
\Omega((\gn_{\pm})_{\lambda}^{\ast}) \in (\gn_{\pm})_{\lambda}^{\ast}\otimes
(\gn_{\pm})_{\lambda}$,
$\Omega^0:=(\Omega(\gh)+\Omega(\gh^{\ast})) \in \gh\otimes\gh^{\ast}+ 
\gh^{\ast}\otimes\gh$. Set
$$ \Omega_{\pm}=\sum_{\lambda}\Omega_{\lambda,\pm}^- + \Omega^0 +
\sum_{\lambda}\Omega_{\lambda,\pm}^+ \in D(\gb_{\pm})\otimes D(\gb_{\pm}).$$
Let $\Omega_{+,ij}$ be the operator on $M$ (or $M^{\ast}$) acting as
$\Omega_+$ on $M(\Lambda_i)\otimes M(\Lambda_j)$, 
($M(\Lambda_i)^{\ast}\otimes M(\Lambda_j)^{\ast}$ respectively) and as the identity
on the other factors. The action of $D(\gb_{\pm})$ on $M(\Lambda_j)$ is the action decribed 
in the previous chapter when $n=1$.
 The KZ equations with additional parameter $\mu\in\gh$, for a function $u(\mu,z)$
 on $n$  variables $z=(z_1,\ldots,z_n)$ taking values in $M$ (or $M^{\ast}$) are
\begin{equation}\label{genKZ}
 \kappa\frac{\partial u}{\partial z_i}=\mu^{(i)}u + \sum_{i\ne j}
 \frac{\Omega_{+,ij}}{z_i - z_j}u, \qquad\qquad i=1,\ldots,n, 
\qquad  \kappa\in\cit.
\end{equation}

Let $\alpha$ be a positive root for $\gg$ and $(y_i^{(\alpha)})$  a basis 
of  $(\gnp)_{\alpha}$. Set $x_i^{(\alpha)}=\tau(y_i^{(\alpha)})$.
Then $((y_i^{(\alpha)})^{\ast})$,
$(x_i^{(\alpha)})$ and $((x_i^{(\alpha)})^{\ast})$ are
 bases of $(\gnp^{\ast})_{\alpha}$, $(\gnm)_{\alpha}$,
 $(\gnm^{\ast})_{\alpha}$, respectively.  Define  operators
 $\Delta_{\pm,\alpha}$ on $M$ (or $M^{\ast}$)
via the formulae $\Delta_{+,\alpha} = \sum_i y_i^{(\alpha)}(y_i^{(\alpha)})^{\ast}$, 
$\Delta_{-,\alpha} = \sum_i (x_i^{(\alpha)})^{\ast} x_i^{(\alpha)}$.
The dynamical differential equations for the function  $u(\mu,z)$ 
with values in $M$ ($M^{\ast}$ 
respectively) are
\begin{equation}\label{DynKZ}
 \kappa\partial_{\mu'}u= \sum_{i= 1}^{n}z_i(\mu')^{(i)}u + 
\sum_{\alpha>0 }\frac{\la\alpha,\mu'\ra}{\la\alpha,\mu\ra}\Delta_{+,\alpha}u,
\qquad\qquad \mu' \in \gh, \qquad \kappa \in \cit.
\end{equation}

\subsection{Properties of the operators  $\Delta_{+,\alpha}$.}
 The properties of the operator $\Omega_{+,ij}$ are thoroughly described in 
\cite{V1}. Now we are going to study the operators  $\Delta_{+,\alpha}$.
\begin{lemma}\label{delta+-} The following diagram is commutative:
\[ \begin{CD}
 M @>{\Delta_{+,\alpha}}>> M\\
 @VV{S}V          @VV{S}V \\
 M^{\ast} @>{-\Delta_{-,\alpha}}>> M^{\ast}
\end{CD} \]
 In particular, the operators  $\Delta_{+,\alpha}$ preserve the kernel of the map 
$S:M\rightarrow M^{\ast}$.
\end{lemma}
{\it{Proof.} }We fix $\alpha$ throughout  this proof and will drop it from the 
notation of the bases. Fix a basis $(u_k)$ of $\ker S: (\gnp)_{\alpha}
\rightarrow (\gnp^{\ast})_{\alpha}$. Complete it to a basis of $(\gnp)_{\alpha}$
by  vectors $(v_l)$. Let $(u_k^{\ast})$, $(v_l^{\ast})$ be the dual basis of 
$(\gnp^{\ast})_{\alpha}$. Moreover, 
\begin{equation}\label{S-matrix}
(v_l^{\ast})=\sum_p (A^{-1})_{lp}S(v_p, . ),\quad
\mbox{ and } (\tau(v_l)^{\ast})=\sum_p (A^{-1})_{lp}S(\tau(v_p), . ),
\end{equation}
where $A=(a_{lp})$ is a nondegenerate matrix with entries $a_{lp}=S(v_l,v_p)$.

 For  $y\in(\gnp)_{\alpha}$, consider the map $y:M\rightarrow M$ via
 the action of  $D(\gbp)$. Let $y\cdot p$ denote the $D(\gbp)$ action, and
 $yp$ denote the standard action of $\gbp$. Then we have
\begin{equation}\label{db+act}
 S(y\cdot p,q)=S(\tau(y)p,q)=S(p,-yq) \qquad\mbox{ for any } p,q\in M.
\end{equation}
 Consider the map $(v_l)^{\ast}:M\rightarrow M$ via the action of $D(\gbp)$. For any
$p,q\in M$ we have
\begin{eqnarray} 
S((v_l)^{\ast}p,q) &=& \la(v_l)^{\ast}p,S(q,\, .\,)\ra = 
     \la(\nu_M^{\ast}(\tau(v_l)^{\ast},.))^{\ast}p,S(q,\, .\,)\ra\label{db+st}\\
& = &\la p, \nu_M^{\ast}(\tau(v_l)^{\ast},S(q,\, .\,))\ra
   =  \la p, \nu_M^{\ast}(\sum_j (A^{-1})_{lj}S(\tau(v_j), . ) ,S(q,\, .\,))\ra\notag\\
& = & \sum_j (A^{-1})_{lj} \la p, S(\tau(v_j)q, .)\ra = 
    \sum_j (A^{-1})_{lj} S(p,\tau(v_j)q).\notag
\end{eqnarray}  
The first three equalities come from the definition of the action of $D(\gbp)$ on $M$, 
the last two from formula~(\ref{S-matrix}). We combine (\ref{db+act}), (\ref{db+st}), and
Corollary~\ref{CDcor} to obtain:
\begin{eqnarray}
S(\Delta_{+,\alpha}p,q)& = &S((\sum_k u_k (u_k^{\ast}) + \sum_l v_l
(v_l^{\ast}))p,q)
=S(\sum_l v_l (v_l^{\ast})p,q)\notag\\
&=& S((v_l^{\ast})p,(-v_l)q)=
-\sum_{j,l} (A^{-1})_{lj} S(p,\tau(v_j)v_lq).\label{db+up}
\end{eqnarray}

 Now we trace the arrows in the alternative direction. For $x\in(\gnm)_{\alpha}$,
consider the map $x:M^{\ast}\rightarrow M^{\ast}$ via the action of  $D(\gbm)$.
Denote this action by '$\cdot$'. Let $xq$ denote the standard action of $\gbm$ on $M$.
We have
\begin{equation}\label{db-act}
\la x\cdot S(p, . ),q\ra = \la S(p, -\tau(\tau(x))\, .\, ),q\ra 
= S(p, -x q) \qquad \mbox{ for all } x\in\gbm;\, p,q\in M.
\end{equation}
 Consider the
map $(\tau(v_l))^{\ast}:M^{\ast}\rightarrow M^{\ast}$ via the action of $D(\gbm)$.
\begin{eqnarray} 
\la(\tau(v_l))^{\ast}S(p,\, . \,),q\ra &=&
\la(\sum_j (A^{-1})_{lj}S(\tau(v_j), . ))S(p,\, .\,),q\ra\label{db-st}\\
& = &\sum_{j} (A^{-1})_{lj}\la S(\tau(v_j)p,\, .\,),q\ra = 
\sum_{j} (A^{-1})_{lj}S(\tau(v_j)p,q).\notag
\end{eqnarray} 
Finally combine (\ref{db-act}), (\ref{db-st}), and Corollary~\ref{CDcor} with $u_k\in\ker S$
to get
\begin{eqnarray}
\lefteqn{\la -\Delta_{-,\alpha}S(p,\,.\,),q\ra =
\la-\left(\sum_k\tau(u_k)^{\ast}\tau(u_k) 
+ \sum_l\tau(v_l)^{\ast}\tau(v_l)\right)S(p,\,.\,),q\ra}\notag\\
&=&-\la\sum_k\tau(u_k)^{\ast}S(p,-\tau(u_k)\,.\,),q\ra - 
\la\sum_l\tau(v_l)^{\ast}S(p,-\tau(v_l)\,.\,),q\ra\notag\\
&=& \la-\sum_k\tau(u_k)^{\ast}S(-\tau(-\tau(u_k))p,\,.\,),q\ra+
\la-\sum_l\tau(v_l)^{\ast}S(-\tau(-\tau(v_l))p,\,.\,),q\ra\notag\\
&=& 0 -\la\sum_{j,l} (A^{-1})_{lj}S(\tau(v_j)(v_l)p,\,.\,),q\ra
=-\sum_{j,l} (A^{-1})_{lj}S(\tau(v_j)(v_l)p,q)\notag\\
&=&-\sum_{j,l} (A^{-1})_{lj}S(p,\tau(v_l)v_jq).\label{db-down}
\end{eqnarray}
Since the matrix $A$ is symmetric (\ref{db+up}) and (\ref{db-down}) prove 
that the diagram is commutative.\qed 

    As a corollary of the Lemma we have that $\Delta_{+,\alpha}$
naturally acts on
$L=M/\ker(S:M\rightarrow M^{\ast})$. We describe this action.
Consider the Kac-Moody algebra $\bar{\gg}=\gg/\ker(S:\gg\rightarrow\gg^{\ast})$.
Let $x\mapsto \bar x$ denote the canonical projections $M\to L$, $\gg\to\bar{\gg}$.
$\ker S$ is an ideal and the form $S$ induces a non-degenerate Killing form on 
$\bar{\gg}$ via the formula $K(x,y)=-S(\tau(x),y)$, see \cite{V1}. $K$ induces a
non-degenerate pairing between root spaces $\bar{\gg}_{\alpha}$ and 
$\bar{\gg}_{-\alpha}$. Let $(e_l^{(\alpha)})$ be a basis of $\bar{\gg}_{\alpha}$,
and let $(f_l^{(\alpha)})$ be the dual basis of $\bar{\gg}_{-\alpha}$ with respect 
to $K$.  Let $\Bar{\Delta}_{\alpha}= \sum_l f_l^{(\alpha)} e_l^{(\alpha)}$.
\begin{corollary}\label{corCD} The following diagram is commutative
\[ \begin{CD}
  M @>{\Delta_{+,\alpha}}>{D(\gbp)-action}>  M\\
  @VVV   @VVV\\
  L @>{\Bar{\Delta}_{\alpha}}>{standard}> L.
\end{CD} \]
\end{corollary}
{\it{ Proof.}} $L(\Lambda)\cong \img \{S:M\rightarrow M^*\}$ via 
$\bar{x}\mapsto S(x,\,.\,)$.
We keep the notation from the Lemma above. Set $w_l=-\sum_j(A^{-1})_{lj}\tau(v_j)$.
From the computation in the Lemma we have
\begin{eqnarray*}
 \overline{\Delta_{+,\alpha}p}&=& S(\Delta_{+,\alpha}p,\,.\,)
=S(\sum_{l}(-\sum_j(A^{-1})_{lj}\tau(v_j))v_lp,\,.\,)\\
&=&S(\sum_lw_lv_lp,\,.\,)=\overline{\sum_{l}w_lv_lp}=\sum_l\bar{w_l}\bar{v_l}\bar{p}
\end{eqnarray*}
Finally notice that the set $(\bar{v_l})$ forms a basis of $\bar{\gg}_{\alpha}$, and
the set $(\bar{w_l})$ forms the dual basis of $\bar{\gg}_{-\alpha}$ with respect to $K$.\qed
\begin{corollary}\label{m2mstar} Fix $\lambda\in\nit^r$.
 Let $m\in M_{\lambda}$, and let $(m_j)$ be a basis of $M_{\lambda}$, and let 
$(m_j^{\ast})$ be the dual basis of $M^{\ast}_{\lambda}$. Then the following 
decomposition holds
$$\Delta_{+,\alpha}m=\sum_j \la-\Delta_{-,\alpha}m_j^{\ast},m\ra m_j.$$
\end{corollary}
{\it{Proof.}} Let $y\in \gbp$ and $x=\tau(y)\in \gbm$. Let $p\in M$, and
$\phi\in M^{\ast}$.
 As in the proof of Lemma~\ref{delta+-} $\la y^{\ast}p,\phi\ra=\la p,x^{\ast}\phi\ra$, where 
 $D(\gbp)$ acts on $M$ and $D(\gbm)$ acts on $M^{\ast}$. Moreover
$\la y\cdot p,\phi\ra= \la p,-x\cdot\phi\ra$, where  $D(\gbp)$ acts on $M$ and 
$D(\gbm)$ acts on $M^{\ast}$. Finally noting that 
$\Delta_{+,\alpha}=\sum_iy_i^{(\alpha)}(y_i^{(\alpha)})^{\ast}$ and
$\Delta_{-,\alpha}=\sum_i\tau(y_i^{(\alpha)})^{\ast}\tau(y_i^{(\alpha)})$
we have
$$\Delta_{+,\alpha}m=\sum_j \la m_j^{\ast},\Delta_{+,\alpha}m\ra m_j=
\sum_j \la-\Delta_{-,\alpha}m_j^{\ast},m\ra m_j.\mbox{\qed}$$

\subsection{An integral form of the dynamical differential equations}
  Our aim now is to rewrite the Dynamical  equations in a form related to the 
hypergeometric solutions.  Fix $\lambda=(m_1,\cdots,m_r)\in\nit^r$. Let $M$ be a tensor 
product of highest weight modules of the Kac-Moody Lie algebra without Serre's relations $\gg$.
 Set $\Lambda=\sum\Lambda_j$, the sum of the respective highest weights. 
Consider the weight space $M_{\lambda}$ of $M$ with weight
$\Lambda - \alpha(\lambda)$, where $\alpha(\lambda)=\sum m_i\alpha_i$. 
Fix a highest-weight vector $v_j$ for each module $M(\Lambda_j)$, 
$j=1,\ldots,n$. To every $I\in P(\lambda,n)$ we associate a vector 
$f_Iv=f_{i_1^1}\cdots f_{i_{s_1}^1}v_1\otimes\cdots\otimes f_{i_1^n}\cdots 
f_{i_{s_n}^n}v_n$ in $M_{\lambda}$, cf. Section~\ref{HG-sol}. 
Note that the vectors $(f_Iv)_{I\in P(\lambda,n)}$ form a basis of $M_{\lambda}$.

 $\gnm$ acts on $M_{\lambda}^{\ast}$ via the $D(\gbm)$ action. Therefore $U(\gnm)$ acts
on  $M_{\lambda}^{\ast}$. Explicitly, $x\cdot \phi(\, .\,)=\phi(-x\, .\,)$ for $x\in\gnm$, 
$\phi\in M^{\ast}$ (cf. Section~\ref{defKZ}), where the action on the left
hand side is 
the $D(\gbm)$ one and the action on the right hand side is the standard one. 
Let $V$ be a vector space freely generated by $f_1,\ldots,f_r$. Therefore we have 
an inclusion of tensor algebras $T(V)\subset T(\gnm)$. Moreover $T(V)$ is an associative
enveloping algebra of the Lie algebra $\gnm$. Since $T(V)$ is a free associative algebra,
$T(V)$ is isomorphic to the universal enveloping algebra $U(\gnm)$. From now on we
will refer to the monomial basis of $T(V)$ as to the monomial basis of $U(\gnm)$, and
to the dual of the monomial basis of $T(V)$ as to the monomial basis of $U(\gnm)^*$.
Rewrite a commutator $x\in\gnm$ as  an element of $U(\gnm)$ in the form
$x=\sum a_jx_j$ 
where $a_j\in\zit$, and $x_j$'s are elements of the monomial basis of $U(\gnm)$. Thus
\begin{equation}\label{gnmact}
 \la x\cdot \phi,\, .\,\ra=\la\phi,-\sum a_jx_j\, .\,\ra\end{equation}

Denote $i:\gnm\rightarrow U(\gnm)$ the inclusion monomorphism. 
Let $\sigma_j\in \{ 1,\ldots,r\}$ for $j=1,\ldots,k$. Let the positive root of $\gg$, 
$\sum_{j=1}^n\alpha_{\sigma_j}$, correspond to the $r$-tuple $\lambda'\in\nit^r$,
i.e. $\alpha(\lambda')=\sum_{j=1}^n\alpha_{\sigma_j}$. Define an element
$\Delta_{\sigma_1,\ldots,\sigma_k}$ of
$(\gnm)_{\lambda'}^*$ via the rule
\begin{equation}\label{deltadef}
\la\Delta_{\sigma_1,\ldots,\sigma_k},x\ra=\la (f_{\sigma_1}\cdots f_{\sigma_k})^{\ast},i(x)\ra 
\mbox{ where } x\in\gnm, \mbox{ and } 
(f_{\sigma_1}\cdots f_{\sigma_k})^{\ast}\in U(\gnm)^{\ast}.
\end{equation}
Thus $ \la\Delta_{\sigma_1,\ldots,\sigma_k},x\ra$ is the coefficient of 
$f_{\sigma_1}\cdots f_{\sigma_k}$ in the decomposition of $i(x)$ into a sum of monomials .

\begin{lemma}\label{deltalem}
Let $\alpha=\sum_{k=1}^r m'_i\alpha_i$ be a positive root of $\gg$. Let $I\in P(\lambda,n)$. Set 
$\lambda'=(m_1',\ldots,m_r')\in\nit^r$. Then we have
\begin{equation}
- \Delta_{-,\alpha} (f_Iv)^{\ast}= \left(\sum_{(i_1,\ldots,i_{m'})\in
P(\lambda',1)}
\Delta_{i_1,\ldots,i_{m'}}(e_{i_{m'}}\cdots e_{i_1}) \right)(f_Iv)^{\ast},
\end{equation} where $\Delta_{i_1,\ldots,i_{m'}}$ acts according to the $D(\gbm)$ action on
$M^{\ast}_{\lambda}$, the product of $e$'s acts on one tensor factor at a time
$(e_{i_{m'}}\cdots e_{i_1})=\sum_{j=1}^n e_{i_{m'}}^{(j)}\cdots e_{i_1}^{(j)}$, and
$e_j$ acts via the standard action (\ref{*action}) on each tensor factor. 
\end{lemma}
{\it Proof.}
 Note that $\alpha=\alpha(\lambda')$. Let
$x\in(\gnm)_{\lambda'}$. If there exists $i$ such that $m'_i>m_i$, then formula~(\ref{gnmact})
implies $x\cdot (f_Iv)^{\ast}=0$ for every $I\in P(\lambda,n)$ because each monomial 
$x_j$ in the $U(\gnm)$ expansion of $x$ has more $f_i$'s than $f_Iv$.

 Now let $m'_i\leq m_i$ for any $1\leq i\leq r$. 
Set $m'=\sum m'_i$. First consider the case $n=1$, $I=(i_1,\ldots,i_m)\in P(\lambda,1)$.
Let $(x_j^{\ast})$ be a basis
of $(\gnm)_{\lambda'}^{\ast}$ such that $x_1^*=\Delta_{i_1,\ldots,i_{m'}}$, and let
$(x_j)$ be the dual basis of $(\gnm)_{\lambda'}$. Formula~(\ref{deltadef}) implies that
the coefficient  of the monomial $f_{i_1}\cdots f_{i_{m'}}$ in the $U(\gnm)$ expansion
of $x_1$ is $1$, and the coefficient of $f_{i_1}\cdots f_{i_{m'}}$ in the  $U(\gnm)$ expansion of
$x_j$ is $0$ for $j>1$. Now use (\ref{gnmact}) to obtain
\begin{align}
 -x_1\cdot(f_Iv)^{\ast}&=(f_{i_{m'+1}}\cdots f_{i_m}v)^{\ast}=
e_{i_{m'}}\cdots e_{i_1} (f_Iv)^{\ast}\notag\\
 -x_j\cdot(f_Iv)^{\ast}&=0 \qquad \mbox{ for } j\geq 1. 
\end{align}
For any element $I'=(i'_1,\ldots,i'_{m'})\in P(\lambda',1)$, such that 
$I'\ne(i_1,\ldots,i_{m'})$ we have \\ $e_{i'_{m'}}\cdots e_{i'_1}(f_Iv)^* =0$.
The proof for $n=1$ is finished.

 Let $n$ be arbitrary natural number. 
$(f_Iv)^*=(f_{i_1^1}\cdots f_{i_{s_1}^1}v_1)^*\otimes\cdots\otimes (f_{i_1^n}\cdots 
f_{i_{s_n}^n}v_n)^*$. Formula (\ref{gnmact}) implies $\la - x\cdot (f_Iv)^*,\, .\,\ra = 
\la (f_Iv)^*, x\,.\,\ra=\la (f_Iv)^*, \sum_{j=1}^nx^{(j)}\,.\,\ra$. 
This and the computation for $n=1$ give 
$$- \Delta_{-,\alpha} (f_Iv)^* = \sum_{(i_1,\ldots,i_{m'})\in P(\lambda',1)}
\sum_{j=1}^n \Delta_{i_1,\ldots,i_{m'}}e_{i_{m'}}^{(j)}\cdots e_{i_1}^{(j)}(f_Iv)^*.\qquad\square$$

 For every positive root $\alpha=\sum m_i^{\alpha}\alpha_i$ of $\gg$, set
$\lambda_{\alpha} = (m_1^{\alpha},\ldots,m_r^{\alpha})\in\nit^r$. 
Now we combine Corollary~\ref{m2mstar} and Lemma~\ref{deltalem} to obtain
 the following form of  the Dynamical KZ equation.
\begin{lemma} 
Let $u(\mu,z)=\sum_{I\in P(\lambda,n)}u_If_Iv$, and let $\mu'\in \gh$ be  a direction 
of differentiation. The Dynamical differential equation (\ref{DynKZ}) is 
equivalent to the equation
\begin{eqnarray}\label{DynKZ*}
\lefteqn{\kappa\partial_{\mu'}u= \sum_{i= 1}^{n}z_i(\mu')^{(i)}u +} \\ 
&+&\sum_{I\in P(\lambda,n)}\sum_{\alpha>0 }\frac{\la\alpha,\mu'\ra}{\la\alpha,\mu\ra}
\sum_{J\in P(\lambda,n)}\la(\sum_{(i_1,\ldots,i_{m'})\in P(\lambda_{\alpha},1)}
\Delta_{i_1,\ldots,i_{m'}}(e_{i_{m'}}\cdots e_{i_1}))(f_Jv)^{\ast},f_Iv\ra u_If_Jv.\notag
\quad\square\end{eqnarray}
\end{lemma}

 \subsection{A symmetrization procedure}\label{symsection}
 The definition of the hypergeometric differential form involves a symmetrization 
procedure, see Section~\ref{HG-sol}. Now we will study the behavior of the operator
$\Delta_{i_1,\ldots,i_{m'}}e_{i_{m'}}\cdots e_{i_1}$ for
$(i_i,\ldots,i_{m'})\in P(\lambda',1)$, where $\alpha=\alpha(\lambda')$,
under the same type of symmetrization procedure.

{\bf Complexes \cite{SV1}} 
For a Lie algebra $\gg$ and a $\gg$-module M, denote by
$C_{\bullet}(\gg,M)$ the standard 
chain complex of $\gg$ with coefficients in $M$. $C_p(\gg,M)=\Lambda^p\gg\otimes M$ and
\begin{align}
 d: g_p\wedge\cdots\wedge g_1\otimes x &= \sum_{i=1}^p(-1)^{i-1}g_p\wedge\cdots\wedge
\widehat{g_i}\wedge\cdots\wedge g_1\otimes g_ix + \notag\\
&+\sum_{1\leq i<j\leq p}(-1)^{i+j}g_p\wedge\cdots\wedge \widehat{g_j}\wedge\cdots
\wedge
\widehat{g_i} \wedge\cdots\wedge g_1\wedge[g_j,g_i]\otimes x.\notag
\end{align}
  Let $\Lambda_1,\ldots,\Lambda_n \in\gh^{\ast}$. Set 
$M=M(\Lambda_1)\otimes\cdots\otimes M(\Lambda_n)$. Consider the
complex $C_{\bullet}(\gnm,M)$. We have the weight decomposition
$C_{\bullet}(\gnm,M)=\oplus_{\lambda\in\nit^r}C_{\bullet}(\gnm,M)_{\lambda}$.

In Section~\ref{defKZ} we recalled a Lie algebra structure on
$\gnm^{\ast}$ and a $\gnm^{\ast}$--module structure on $M^{\ast}$.
Let $C_{\bullet}(\gnm^{\ast},M^{\ast})$ be the corresponding standard chain
complex: $C_{\bullet}(\gnm^{\ast},M^{\ast})=\oplus_{\lambda\in\nit^r}
C_{\bullet}(\gnm^{\ast},M^{\ast})_{\lambda}$. The covariant form induces a
graded homomorphism of complexes,
$S:C_{\bullet}(\gnm,M)\rightarrow C_{\bullet}(\gnm^{\ast},M^{\ast})$, see \cite{SV1}.

Let $\lambda=(m_1,\ldots,m_r)\in\nit^r$, and $m=\sum m_i$. Define a subgroup 
$\Sigma_{\lambda}$ of the symmetric group $\Sigma_m$ via the direct product
 $\Sigma_{\lambda}=\Sigma_{m_1}\times\cdots\times\Sigma_{m_r}$, where $\Sigma_{m_j}$ 
permutes the set of indices 
$\{\sum_{p=1}^{j-1}m_p+1,\ldots,\sum_{p=1}^{j-1}m_p+m_j\}$.
Introduce a free Lie algebra
$\tnm$ on generators $\tilde{f_1},\ldots,\tilde{f_m}$. Define a map of Lie algebras
$\gnm\rightarrow\tnm$ by setting $f_i\mapsto \sum_{j=1}^{m_i} \tilde{f}_{m(i)+j}$,
where $m(i)=m_1+\cdots+m_{i-1}$. It induces a map of complexes
$C_{\bullet}(\gnm,U(\gnm)^{\otimes n})\rightarrow 
C_{\bullet}(\tnm,U(\tnm)^{\otimes n})^{\Sigma_{\lambda}}$. Set 
$ \widetilde{\lambda}=(\underbrace{1,1,\ldots,1}_{m})$. Let
\begin{equation}
 s: C_{\bullet}(\gnm,U(\gnm)^{\otimes n})_{\lambda}\rightarrow 
C_{\bullet}(\tnm,U(\tnm)^{\otimes n})^{\Sigma_{\lambda}}_{\widetilde{\lambda}}
\end{equation}
be the previous map composed with the projection on the $\widetilde{\lambda}$--
component.

 On the other hand there is a map of Lie algebras $\tnm\rightarrow\gnm$ defined
by $\widetilde{f}_j\mapsto f_i$, for $m(i)<j\leq m(i+1)$. It induces the map
$\pi_{\lambda}:C_{\bullet}(\tnm,U(\tnm)^{\otimes n})_{\widetilde{\lambda}}
\rightarrow C_{\bullet}(\gnm,U(\gnm)^{\otimes n})_{\lambda}$.
Note that $s(y)$ equals the sum over the preimages  of $y$ under $\pi_{\lambda}$, for
any $y\in C_{\bullet}(\gnm,U(\gnm)^{\otimes n})_{\lambda}$. Each such preimage is 
uniquely described by an element $\sigma\in\Sigma_{\lambda}$.

{\it Example.} Let $n=1$ and $I=(i_1,\ldots,i_m)\in P(\lambda,1)$. Consider
$f_I$ as an element of $C_0(\gnm,U(\gnm))_{\lambda}$. Then 
$s(f_I)=\sum_{\sigma}\tf_{\sigma_1}\cdots\tf_{\sigma_m}$, where the sum is over
the set $\{\sigma\in \S(I)\}\cong\Sigma_{\lambda}$.
\begin{lemma}\label{l-symf} Let $n=1$, $I\in P(\lambda,1)$, and $m'\leq m$. Then
the map $s^*:(C_1(\tnm,U(\tnm))_{\widetilde{\lambda}}^{\Sigma_{\lambda}})^*\to
(C_1(\gnm,U(\gnm))_{\lambda})^*$ has the following property
\begin{equation}\label{symf}
s^{\ast}\left(\frac{1}{|\Sigma_{\lambda}|}\sum_{\sigma\in \S(I)}
\widetilde{\Delta}_{\sigma_{1},\ldots,\sigma_{m'}}\otimes
(\widetilde{f}_{\sigma_{m'+1}}\cdots\widetilde{f}_{\sigma_m})^{\ast}\right)=
\Delta_{i_1,\ldots,i_{m'}}\otimes (f_{i_{m'+1}}\cdots f_{i_m})^{\ast}
\end{equation}
\end{lemma}
{\bf Proof.} Choose a basis $(x_j^{\ast})$ of $\gnm^{\ast}$ such that 
$x_1^{\ast}=\Delta_{i_1,\ldots,i_{m'}}$ and let $(x_j)$ be the dual basis of $\gnm$.
Take a tensor product of the basis $(x_j)$ with the monomial basis in $U(\gnm)$ to 
get a basis in $\gnm\otimes U(\gnm)$. Let us compare the two sides of (\ref{symf})
on that basis. 
$\Delta_{i_1,\ldots,i_{m'}}\otimes (f_{i_{m'+1}}\cdots f_{i_m})^{\ast}(x_1\otimes
f_{i_{m'+1}}\cdots f_{i_m})=1$. The right hand side is zero on any other element 
of the basis.

 Let $i:\gnm\otimes U(\gnm)\rightarrow U(\gnm)\otimes U(\gnm)$, and
$\widetilde{i}:\tnm\otimes U(\tnm)\rightarrow U(\tnm)\otimes U(\tnm)$ be the natural
inclusion maps. Clearly $\widetilde{i}\circ s=s\circ i$. Moreover,
$\la\widetilde{\Delta}_{\sigma_{i_1},\ldots,\sigma_{i_{m'}}},s(x)\ra=
\la (\widetilde{f}_{\sigma_{i_{1}}}\cdots\widetilde{f}_{\sigma_{i_{m'}}})^{\ast}
,(\widetilde{i}\circ s(x))\ra$ by definition for $x\in \gnm$.
 For a fixed $\sigma\in \Sigma_{\lambda}$ we have
\begin{eqnarray*}
\lefteqn{\la s^{\ast}(\widetilde{\Delta}_{\sigma_{1},\ldots,\sigma_{m'}}\otimes
(\widetilde{f}_{\sigma_{m'+1}}\cdots\widetilde{f}_{\sigma_{m}})^{\ast})
,x_1\otimes f_{i_{m'+1}}\cdots f_{i_m}\ra=}\\
&=&\la\widetilde{\Delta}_{\sigma_{i_1},\ldots,\sigma_{i_{m'}}}\otimes
(\widetilde{f}_{\sigma_{i_{m'+1}}}\cdots\widetilde{f}_{\sigma_{i_m}})^{\ast}
,s(x_1\otimes f_{i_{m'+1}}\cdots f_{i_m})\ra\\
&=&\la(\widetilde{f}_{\sigma_1}\cdots\widetilde{f}_{\sigma_{m'}})^{\ast}
\circ\widetilde{i}\otimes
(\widetilde{f}_{\sigma_{m'+1}}\cdots\widetilde{f}_{\sigma_m})^{\ast}
,s(x_1\otimes  \widetilde{f}_{i_{m'+1}}\cdots \widetilde{f}_{i_m})\ra\\
&=&\la(\widetilde{f}_{\sigma_{1}}\cdots\widetilde{f}_{\sigma_{m'}})^{\ast}
\otimes (\widetilde{f}_{\sigma_{m'+1}}\cdots\widetilde{f}_{\sigma_m})^{\ast}
,s(i(x_1)\otimes  \widetilde{f}_{i_{m'+1}}\cdots \widetilde{f}_{i_m})\ra
\end{eqnarray*}
The duality of the bases implies
$ i(x_1\otimes f_{i_{m'+1}}\cdots f_{i_m})
=(f_{i_1}\cdots f_{i_{m'}}+\mbox{ other monomials})\otimes f_{i_{m'+1}}\cdots f_{i_m}$,
and $s\circ i(x_1\otimes f_{i_{m'+1}}\cdots f_{i_m})=
(\widetilde{f}_{\sigma_1}\cdots\widetilde{f}_{\sigma_{m'}}+\mbox{ other monomials })
\otimes\tf_{\sigma_{m'+1}}\cdots \tf_{\sigma_m}$. Therefore,
$\la s^{\ast}(\widetilde{\Delta}_{\sigma_{1},\ldots,\sigma_{m'}}\otimes
 (\widetilde{f}_{\sigma_{m'+1}}\cdots\widetilde{f}_{\sigma_{m}})^{\ast})
 ,x_1\otimes f_{i_{m'+1}}\cdots f_{i_m}\ra=1$. The same way we check that the left hand 
side is zero on the other basis elements.  \qed
\begin{corollary}\label{symcor}
 Let $\pi$ be the restriction of the projection $\pi_{\lambda}$ to the 
subspace $C_.(\tnm,U(\tnm)^{\otimes n})^{\Sigma_{\lambda}}_{\widetilde{\lambda}}$ of
$C_.(\tnm,U(\tnm)^{\otimes n})_{\widetilde{\lambda}}$. 
Let $J=(j_1,\ldots,j_m)\in P(\lambda,n)$ and
$I=(i_1,\ldots,i_{m'})\in P(\lambda_{\alpha},1)$ for a positive root $\alpha$. Then we have
\begin{equation}\label{pi-star}
(\sum_{\tau \in \S(I)}\widetilde{\Delta}_{\tau_1,\ldots,\tau_{m'}}
(\te_{\tau_{m'}}\ldots\te_{\tau_1}))
(\sum_{\sigma \in \S(J)}(\tf_{\sigma_1}\cdots
\tf_{\sigma_{m}})^{\ast})=\pi^{\ast}((\Delta_{i_1,\ldots,i_{m'}}(e_{i_{m'}}\ldots 
e_{i_1}))(f_{j_1}\cdots f_{j_m})^{\ast}).
\end{equation}
\end{corollary}
{\it Proof.} Assume $n=1$. The general case follows from this one because the operators
$(\te_{\tau_{m'}}\ldots\te_{\tau_1})$ and $(e_{i_{m'}}\ldots e_{i_1})$
 act on one tensor factor at a time. Since $s\circ \pi = |\Sigma_{\lambda}|\mathbf{id}$ 
we have $\pi^{\ast}\circ s^{\ast}=|\Sigma_{\lambda}|\mathbf{id}$. $s^{\ast}$ and 
$\pi^{\ast}$ are maps of complexes, i.e. they commute with the corresponding 
differentials. Thus, Lemma~\ref{l-symf} implies after applying differentials and 
taking $\pi^{\ast}$ from both sides
$$\pi^{\ast}(\Delta_{j_1,\ldots,j_{m'}}(f_{j_{m'+1}}\cdots f_{j_m})^{\ast})=
\sum_{\sigma \in \S(J)}\widetilde{\Delta}_{\sigma_1,\ldots,\sigma_{m'}}
(\tf_{\sigma_{m'+1}}\cdots\tf_{\sigma_{m}})^{\ast}.$$
The right hand side of our formula is non-zero if and only if
$m'\leq m$ and $i_k=j_k$ for $1\leq k \leq m'$. Thus we compute
\begin{align}
& \pi^{\ast}(\Delta_{i_1,\ldots,i_{m'}}e_{i_{m'}}\ldots 
e_{i_1}(f_{j_1}\cdots f_{j_m})^{\ast})=
\pi^{\ast}(\Delta_{j_1,\ldots,j_{m'}}(f_{j_{m'+1}}\cdots f_{j_m})^{\ast})\notag\\
&= \sum_{\sigma \in \S(J)}\widetilde{\Delta}_{\sigma_1,\ldots,\sigma_{m'}}
(\tf_{\sigma_{m'+1}}\cdots\tf_{\sigma_{m}})^{\ast}=
\sum_{\sigma \in \S(J)}\widetilde{\Delta}_{\sigma_1,\ldots,\sigma_{m'}}
\te_{\sigma_{m'}}\cdots\te_{\sigma_1}
(\tf_{\sigma_{1}}\cdots\tf_{\sigma_{m}})^{\ast}\notag
\end{align}
Note that $(\sigma_1,\ldots,\sigma_{m'})\in \S(I)$. For any other $\tau \in \S(I)$,
we have $\te_{\tau_{m'}}\cdots\te_{\tau_1}
(\tf_{\sigma_{1}}\cdots\tf_{\sigma_{m}})^{\ast}=0$. Therefore we rewrite the last equality 
in the form (\ref{pi-star}).\qed

Our last step is to show that the Dynamical equations for any Lie algebra in any 
weight space follow from the Dynamical equations for any Lie algebra in a weight space 
with weight  $\widetilde{\lambda}=(1,1,\ldots,1)$.

Fix a finite dimensional complex vector space $\gh$, a non-degenerate symmetric bilinear form
$(\,.\,,\,.\,)$ on $\gh$, and a set of linearly independent ``simple roots'' 
$\alpha_1,\ldots,\alpha_r \in\gh^{\ast}$. Consider the corresponding 
Kac-Moody Lie algebra without Serre's  relations, $\gg$, defined at the beginning 
of Section~\ref{defKZ}. Recall that $b:\gh\rightarrow\gh^{\ast}$ denotes
the isomorphism
induced by the bilinear form, $h_i=b^{-1}(\alpha_i)$, and and the form $(\,.\,,\,.\,)$ is
transfered to $\gh^{\ast}$ via the map $b$.
 
 Fix $\lambda=(m_1,\ldots,m_r)\in\nit^r$. Consider the corresponding positive root 
$\alpha(\lambda)=\sum_{i=1}^rm_i\alpha_i$ of $\gg$. Up to reordering of the 
$\alpha$'s we can assume that $\alpha(\lambda)=\sum_{i=1}^{p}m_i\alpha_i$ where
$m_i>0$ for $1\leq i\leq p$ and $p\leq r$ is fixed.
The corresponding coloring function is
$c_{\lambda}:\{1,\ldots,m=\sum m_i\}\rightarrow\{1,\ldots,p\}$.
 We use the following linear 
algebraic fact when symmetrizing.
\begin{proposition} \label{lin-mono}
Let $\gh$ be a finite dimensional vector space with a non-degenerate symmetric
bilinear form $(\,.\,,\,.\,)$, and a set of linearly independent vectors $(h_i)_{i=1}^r
\subset\gh$.
Then there exists a finite dimensional vector space $\tgh$ with a non-degenerate 
symmetric bilinear form
$(\,.\,,\,.\,)_1$, a set of linearly independent vectors $(\thh_j)_{j=1}^m\subset\tgh$, 
and a monomorphism $s_h:\gh\rightarrow\tgh$ such that\\
(a) $s_h(h_i)=\frac{1}{m_i}\sum_{j=1}^{m_i}\thh_{m(i)+j}$, where $m(i)=m_1+\cdots+m_{i-1}$,
$i=1,\ldots,p$;\\ 
(b) $(\thh_j,s_h(h'))_1=(h_{c(j)},h')$
 and $(h',h'')=(s_h(h'),s_h(h''))_1$ for any $h',h''\in\gh$, $j=1\ldots,m$.
\end{proposition} 
{\it Proof.} Let $q=\dim \gh$. Complete the set $h_1,\ldots,h_r$ to a basis 
$h_1,\ldots,h_r,h_{r+1},\ldots,h_q$ of $\gh$. U
Consider a complex linear space 
$\tgh'=\cit\{\thh_1,\ldots,\thh_m,\thh_{m+1},\ldots,\thh_{m+q-p}\}$. Extend the coloring 
function $c:\{1,\ldots,m+q-p\}\rightarrow\{1,\ldots,q\}$ setting 
$c(m+j)=r+j$ for $j=1,\ldots,q-p$.Define a symmetric
degenerate bilinear form on $\tgh'$ by the rules
$(\thh_j,\thh_k)_1=(h_{c(j)},h_{c(k)})$ for $1\leq j,k\leq m+q-p$. The rank of the form is
$q$ and the dimension of its kernel is $m-p$. There exists an extension $\tgh$ 
of the vector space $\tgh'$ and an extension of  $(\,.\,,\,.\,)_1$ to a non-degenerate 
symmetric bilinear form on $\tgh$. 

Define a monomorphism $s_h$ by $s_h(h_i)=\frac{1}{m_i}\sum_{j=1}^{m_i}\thh_{m(i)+j}$, 
where $m(i)=m_1+\cdots+m_{i-1}$, $i=1,\ldots,q$, and note that $s_h(h_{p+j})=\thh_{m+j}$ for
$j=1,\ldots,q-p$. Now checking (b) on a basis is straightforward.\qed

Set $\talp_j=(\thh_j,\,.\,)_1\in\tgh^{\ast}$ for $j=1\ldots,m$. Consider $\tgg$, 
a Kac-Moody Lie algebra 
without Serre's relations corresponding to the data $\tgh$, $(\,.\,,\,.\,)_1$, and
$(\talp_j)_{j=1}^m$. Note that $1\leq j\leq m$ 
implies $\la\talp_j,s_h(h')\ra=(\thh_j,s_h(h'))_1=(h_{c(j)},h')=\la\alpha_{c(j)},h'\ra$ 
for any $h'\in\gh$.
 
 Let $M=M(\Lambda_1)\otimes\cdots\otimes M(\Lambda_n)$ be a tensor product 
of Verma modules for $\gg$ with corresponding highest weights $\Lambda_1,\ldots,\Lambda_n
\in \gh^{\ast}$. Since $s_h$ is a monomorphism, $s_h^{\ast}:\tgh^{\ast}\rightarrow
\gh^{\ast}$ is a linear epimorphism. Choose highest weights 
$\tLam_1,\ldots,\tLam_n\in\tgh^{\ast}$ such that $s_h^{\ast}(\tLam_j)=\Lambda_j$ for
$1\leq j\leq n$, and consider the corresponding tensor product of Verma modules for
$\tgg$,  $\tM=\tM(\tLam_1)\otimes\cdots\otimes\tM(\tLam_n)$.

\begin{lemma}\label{sym-lemma} Let $\widetilde{\lambda}=(\underbrace{1,1,\ldots,1}_{m})$. 
Let $\tu(\tmu,z)=\sum_{K\in P(\widetilde{\lambda},n)}\tu_K\tf_K$ be a hypergeometric 
solution of the Dynamical equations with values in the $\widetilde{\lambda}$ 
weight space of a $\widetilde{\gg}$-module $\tM\cong U(\tnm)^{\otimes n}$.
 Then $u(\mu,z)=\pi (\tu)(s_h(\mu),z)$ is a hypergeometric solution of the Dynamical equations 
with values in the $\lambda$ weight space of a $\gg$-module $M\cong U(\gnm)^{\otimes n}$,
i.e. $u=\sum_{I\in P(\lambda,n)} u_I f_I$.
\end{lemma}
{\it  Proof.} Note that by definition $P(\tlam,n)= \cup_{I\in P(\lambda,n)}\{K\in \S(I)\}$.
From the definition of the hypergeometric differential form, see Section~\ref{HG-sol},
it follows that $u_I=\sum_{k\in \S(I)}\tu_K$. Therefore
\begin{equation}\label{pi-tu} 
\pi (\tu)=\pi(\sum_{I\in P(\lambda,n)}\sum_{K\in \S(I)} \tu_k \tf_K)=
\sum_{I\in P(\lambda,n)}(\sum_{K\in \S(I)} \tu_k)f_I=\sum_{I\in P(\lambda,n)}u_I f_I.
\end{equation} 
Fix a point $\mu\in\gh$ and a direction of differentiation
$\mu'\in\gh$. Denote $\tmu=s_h(\mu)$ and $\tmu'=s_h(\mu')$. Since 
\begin{align}
&\partial_{\mu'}\exp(-\sum_{j=1}^m \la\alpha_{c(j)},\mu\ra t_j 
+\sum_{l=1}^n\la\Lambda_l,\mu\ra)=\notag\\
&=(-\sum_{j=1}^m \la\alpha_{c(j)},\mu'\ra t_j +\sum_{l=1}^n\la\Lambda_l,\mu'\ra)
\exp(-\sum_{j=1}^m \la\alpha_{c(j)},\mu\ra t_j +\sum_{l=1}^n\la\Lambda_l,\mu\ra)\notag
\end{align}
\begin{align}
&=(-\sum_{j=1}^m \la\talp_j,\tmu'\ra t_j +\sum_{l=1}^n\la\tLam_l,\tmu'\ra)
\exp(-\sum_{j=1}^m \la\talp_j,\tmu\ra t_j +\sum_{l=1}^n\la\tLam_l,\tmu\ra)\notag\\
&=\partial_{\tmu'}\exp(-\sum_{j=1}^m \la\talp_j,\tmu\ra t_j 
+\sum_{l=1}^n\la\tLam_l,\tmu\ra),
\end{align}
we have $\pi(\partial_{\tmu'}\tu(\tmu,z))=\partial_{\mu'}\pi(\tu)(\mu,z)$.
If $I=(I_1,\ldots,I_n)\in P(\lambda,n)$ and $K=(K_1,\ldots,K_n)\in \S(I)$,
where $K_j=(k^j_1,\ldots,k^j_{s_j})$ and $I_j=(i^j_1,\ldots,i^j_{s_j})$, then 
\begin{align}
&\tmu'{}^{(j)}\tf_K=\la\tLam_j-\sum_{l=1}^{s_j}\talp_{k_l^j},\tmu'\ra\tf_K=
\la\Lambda_j-\sum_{l=1}^{s_j}\alpha_{i_l^j},\mu'\ra\tf_K, \notag\\
&\pi(\tmu'{}^{(j)}\sum_{K\in \S(I)}\tu_K \tf_K)=
\la\Lambda_j-\sum_{l=1}^{s_j}\alpha_{i_l^j},\mu'\ra(\sum_{K\in \S(I)}\tu_K\pi(\tf_K))\notag\\
&=\la\Lambda_j-\sum_{l=1}^{s_j}\alpha_{i_l^j},\mu'\ra(\sum_{K\in \S(I)}\tu_K) f_I
=\la\Lambda_j-\sum_{l=1}^{s_j}\alpha_{i_l^j},\mu'\ra u_I f_I=\mu'{}^{(j)} f_I\label{mu-part}
\end{align}
Combine formulae (\ref{pi-tu}) and (\ref{mu-part}) to obtain
\begin{equation}\label{z-part}
\pi(\sum_{j=1}^n z_j\tmu'{}^{(j)}\tu)=\sum_{j=1}^n z_j\mu'{}^{(j)}u.
\end{equation} 

Let $\alpha=\sum_{i=1}^r m'_i\alpha_j$ be a positive root for $\gg$. Lemma~\ref{deltalem}
 gives a  necessary condition for a non-zero action of
$\Delta_{-,\alpha}$ on 
$M_{\lambda}^{\ast}$. Namely $m'_i\leq m_i$ for all $i=1,\ldots,r$. Analogously,
for a positive root $\talp=\sum_{j=1}^m p_j\talp_j$ of $\tgg$ a  necessary 
condition for a non-zero action of $\Delta_{-,\talp}$ on 
$\tM_{\tlam}^{\ast}$ is $p_j=0,1$ for $j=1,\ldots,m$. Call all such $\alpha$'s ( $\talp$'s)
$\lambda$-admissible ($\tlam$-admissible). Since $s_h^{\ast}(\talp_j)=\alpha_{c(j)}$, 
$s_h^{\ast}$ maps the set of $\tlam$-admissible roots of $\tgg$ onto the set of
$\lambda$-admissible roots of $\gg$.

 Let $\alpha=\sum m'_i\alpha_i$ be a $\lambda$-admissible root for $\gg$ and $m'=\sum m'_i$. 
For any $\talp$, such that  $s_h^{\ast}(\talp)=\alpha$, we have 
$\la\talp,\tmu'\ra/\la\talp,\tmu\ra=\la\alpha,\mu'\ra/\la\alpha,\mu\ra$.
Consider $\pi(\sum_{\talp,\,s_h^{\ast}(\talp)=\alpha }
\frac{\la\talp,\tmu'\ra}{\la\talp,\tmu\ra}\Delta_{+,\talp}\tu)$, where the sum is over 
$\tlam$-admissible roots.
 Corollary~\ref{m2mstar} applied to the basis $(\tf_K)_{K\in P(\tlam,n)}$ of $M_{\lambda}$ 
gives
\begin{align}
 &\frac{\la\alpha,\mu\ra}{\la\alpha,\mu'\ra}\pi(\sum_{\talp,\,s_h^{\ast}(\talp)
=\alpha } \frac{\la\talp,\tmu'\ra}{\la\talp,\tmu\ra}
\Delta_{+,\talp}\tu)=\pi(\sum_{K\in P(\tlam,n)}
\la-\sum_{\talp,\,s_h^{\ast}(\talp)=\alpha }\Delta_{-,\talp}(\tf_K)^{\ast},\tu\ra\tf_K)=\notag\\
&=\sum_{I\in P(\lambda,n)}\la(-\sum_{\talp,\,s_h^{\ast}(\talp)=\alpha }
\Delta_{-,\talp})(\sum_{K\in \S(I)}(\tf_K)^{\ast}),\tu\ra f_I\label{pppppp}
\end{align}
 Lemma~\ref{deltalem} asserts that
$$-\sum_{\talp,\,s_h^{\ast}(\talp)=\alpha }
\Delta_{-,\talp}= \sum_{\talp,\,s_h^{\ast}(\talp)=\alpha }
\left(\sum_{(l_1,\ldots,l_{m'})\in P(\tlam_{\talp},1)}
\widetilde{\Delta}_{l_1,\ldots,l_{m'}}(\te_{l_{m'}}\ldots\te_{l_1})\right).$$
Rearrange the summation using that sum over 
$(l_1,\ldots,l_{m'})\in P(\tlam_{\talp},1)$ such that $s_h^{\ast}(\talp)=\alpha$ equals
 the sum over $(p_1,\ldots,p_{m'}) \in \S(J)$ such that $J=(j_1,\ldots,l_{m'})\in 
P(\lambda_{\alpha},1)$. 
Combine such rearrangement with Lemma~\ref{deltalem} and Corollary~\ref{symcor} to simplify
formula~(\ref{pppppp}).
\begin{align}
 &\frac{\la\alpha,\mu\ra}{\la\alpha,\mu'\ra}\pi(\sum_{\talp,\,s_h^{\ast}(\talp)
=\alpha } \frac{\la\talp,\tmu'\ra}{\la\talp,\tmu\ra}\Delta_{+,\talp}\tu)
=\notag\\
&=\sum_{I\in P(\lambda,n)}\la\sum_{J\in P(\lambda_{\alpha},1)}
(\sum_{(p_1,\ldots,p_{m'}) \in \S(J)}
\widetilde{\Delta}_{p_1,\ldots,p_{m'}}(\te_{p_{m'}}\ldots\te_{p_1}))
(\sum_{K\in \S(I)}(\tf_K)^{\ast}),\tu\ra f_I\notag\\
&=\sum_{I\in P(\lambda,n)}\la\sum_{J\in P(\lambda_{\alpha},1)}\pi^{\ast}(
\Delta_{j_1,\ldots,j_{m'}} (e_{j_{m'}}\ldots e_{j_1})(f_I)^{\ast}),\tu\ra f_I\notag\\
&=\sum_{I\in P(\lambda,n)}\la\sum_{J\in P(\lambda_{\alpha},1)}
\Delta_{j_1,\ldots,j_{m'}} (e_{j_{m'}}\ldots e_{j_1})(f_I)^{\ast},\pi(\tu)\ra f_I\notag\\
&=\sum_{I\in P(\lambda,n)}\la-\Delta_{-,\alpha}(f_I)^{\ast},u\ra f_I=
\Delta_{+,\alpha}u.\label{delta-part}
\end{align}
Finally (\ref{z-part}) and (\ref{delta-part}) imply
\begin{equation}
\partial_{\mu'}u=\pi(\partial_{\tmu'}\tu)=\pi((\sum_{j=1}^n z_j\tmu'{}^{(j)}+
\sum_{\talp>0}\frac{\la\talp,\tmu'\ra}{\la\talp,\tmu\ra}\Delta_{+,\talp})\tu)=
(\sum_{j=1}^n z_j\mu'{}^{(j)}+
\sum_{\alpha>0}\frac{\la\alpha,\mu'\ra}{\la\alpha,\mu\ra}\Delta_{+,\alpha})u.
\quad\square
\end{equation}

\subsection{The proof of Lemma~\ref{CDlemma}}\label{proofCDlem}
   Recall that the linear map $\nu_M: M\rightarrow \gbm\otimes M$ has the following 
property $\nu_M(x)=\displaystyle\frac12(b^{-1}(\Lambda-\alpha(\lambda)))\otimes x  +
\nu_{M-}(x)$, where $x\in M_{\lambda}$, $\nu_{M-}(x)\in\gnm\otimes M$,
and
$b^{-1}:\gh^{\ast}\rightarrow\gh$ is defined at the beginning of Section~\ref{defKZ}.
Let $a\in\gh$, $x\in M_{\lambda}$. Since $S(b\otimes x,a\otimes y)=S(a,b)S(x,y)$
 for any $b\in \gg$, $y\in M$,  and $\gh$ is orthogonal to $\gnm$ with respect 
to $S$, and $(\,.\, , \,.\,)$ coincides with $S$ on $\gh$ we have
\begin{equation}
S(\nu_{M-}(x),a\otimes y) =
S(\frac12(b^{-1}(\Lambda-\alpha(\lambda))),a)S(x,y)=
\frac12S(x,ay).
\end{equation}
This proves the first equality in the Lemma. To prove the second part for a monomial
$x$, we use double induction by the number  of tensor factors and the number of $f's$
in $x$. We use $\nu_{M-}(x)$ instead of $\nu_M(x)$ because of the
orthogonality 
mentioned above.

Let $M=M(\Lambda_1)$ be a highest weight module of $\gg$ with a highest vector $v$,
and $a\in\gnm$, $y\in M$. Since $S(\nu_{M-}(v))=0$,  we have
$S(\nu_{M-}(v),a\otimes y)=0=S(v,ay)$. The inductive step is as follows.
Assume $S(\nu_{M-}(x),a\otimes y)=S(x,ay)$. Then
\begin{eqnarray}\label{1fact}
 \lefteqn{ S(\nu_{M-}(f_i x),a\otimes y) = }\notag\\
&=& S(f_i\otimes h_i x,a\otimes y) + S(f_i \nu_{M-}(x), a\otimes y)
 = S(f_i,a)S(h_i x,y) + S(\nu_{M-}(x),e_i(a\otimes y))\notag\\
&=& S(f_i,a)S(h_i x,y) +S(\nu_{M-}(x),[e_i,a]\otimes
y)+S(\nu_{M-}(x),a\otimes e_i y)\notag\\
&=& S(f_i,a)S(h_i x,y) +S(\nu_{M-}(x),[e_i,a]\otimes y)+
S(x,ae_i y)\notag\\
&=& S(f_i,a)S(h_i x,y) +S(\nu_{M-}(x),[e_i,a]\otimes y)+
S(x, e_iay)-S(x,[e_i,a] y)\notag\\
&=& S(f_i,a)S(h_i x,y) +S(\nu_{M-}(x),[e_i,a]\otimes y)+
S(f_ix,ay)-S(x,[e_i,a] y)
\end{eqnarray}
If $a=f_i$ ,then (\ref{1fact}) and the properties $S(x,hy)=S(hx,y)$, 
$S(\nu_{M-}(x),h\otimes y)=0$ for $h\in\gh$ imply
$$ S(\nu_{M-}(f_i x),a\otimes y) = 
  S(h_i x,y) +S(\nu_{M-}(x),h_i\otimes
y)+S(f_ix,ay)-S(x,h_iy)=S(f_ix,y).$$
If $a$ is orthogonal to $f_i$ with respect to $S$, then (\ref{1fact}) and the inductive 
hypothesis give
$$ S(\nu_{M-}(f_i x),a\otimes y) = 0 + S(x,[e_i,a]
y)+S(f_ix,ay)-S(x,[e_i,a] y)
=S(f_ix,a).$$
Thus the statement is proved for one tensor factor.

 Assume that $S(\nu_{M-}(f_i x),a\otimes y)=S(x,ay)$ for a module $M$,
which is a tensor 
product of up to $n-1$ tensor factors $(n\geq 2)$.

 Let $M=M(\Lambda_1)\otimes\cdots\otimes M(\Lambda_n)$. Recall that
$\nu_{M-}(x)=\sum_{k=1}^n\nu^{(k)}_M(x)_-$, where 
$\nu^{(k)}_M(f_i^{(j)}x)_- =f_i^{(j)}\nu^{(k)}_M(x)_-$ for  $ k\ne j$, 
 and $\nu^{(k)}_M(f_i^{(k)}x)_-=f_i\otimes h_i^{(k)}x +
f_i^{(k)}\nu^{(k)}_M(x)_-$. The following commutation relations will be useful.
$S(f_i^{(j)}\nu_{M-}^{(k)}( x),a\otimes y)=S(\nu_{M-}^{(k)}(x),a\otimes
e_i^{(j)}y)$,
 for $j\ne k$, and
$S(f_i^{(j)}\nu^{(j)}_M(x)_-,a\otimes y)=S(\nu_{M-}^{(j)}(x),[e_i,a]\otimes y)+
S(\nu_{M-}^{(j)}( x),a\otimes e_i^{(j)}y)$. Both equalities are
corollaries of the
 Lemma~\ref{CDlemma} for one tensor factor, and the definition of $S$, e.g. 
$S(x_1\otimes\cdots\otimes x_{j-1}\otimes\nu_{M-}(x_j)\otimes x_{j+1} 
\otimes\cdots\otimes x_n, a\otimes y_1\otimes\cdots\otimes y_j \otimes\cdots\otimes y_n)=
S(\nu_{M-}(x_j),a\otimes y_j)\prod_{k\ne j}S(x_k,y_k).$
In all formulae $1\leq i\leq r$, and $1\leq j,k\leq n$, and the upper script indicates
the tensor factor where the action is applied. Let $a\in\gnm$. The base for the induction
is exactly as for $n=1$. 
 The inductive step is as follows.
\begin{eqnarray}\label{mult1}
 \lefteqn{S(\nu_{M-}(f_i^{(j)} x),a\otimes y) 
= \sum_{k\neq j}S(\nu_{M-}^{(k)}(f_i^{(j)} x),a\otimes y)+
S(\nu_{M-}^{(j)}(f_i^{(j)} x),a\otimes y)}\notag\\
&=& \sum_{k\neq j}S(f_i^{(j)}\nu_{M-}^{(k)}( x),a\otimes y)+
S(f_i\otimes h_i^{(j)}x,a\otimes y) +S(f_i^{(j)}\nu^{(j)}_M(x)_-,a\otimes y)\notag\\
&=& S(f_i,a)S(h_i^{(j)}x,y) +S(\nu_{M-}^{(j)}( x),[e_i,a]\otimes y)+
\sum_{k=1}^nS(\nu_{M-}^{(k)}( x),a\otimes e_i^{(j)}y).
\end{eqnarray}
The result for one tensor factor gives
\begin{align}\label{mult2}
\sum_kS(\nu_{M-}^{(k)}( x),a\otimes e_i^{(j)}y) & = \sum_k S( x,a^{(k)}
e_i^{(j)}y)=
\sum_k S( x,e_i^{(j)}a^{(k)}y) -S( x,[e_i,a]^{(j)}y)\notag\\
&= S(f_i^{(j)}x,ay) - S(x,[e_i,a]^{(j)}y).
\end{align}
If $a=f_i$ ,then (\ref{mult1}), (\ref{mult2}) and the properties $S(x,h^{(k)}y)=S(h^{(k)}x,y)$, 
$S(\nu_{M-}^{(k)}(x),h\otimes y)=0$ for $h\in\gh$, $k=1,\ldots,n$ imply
$$S(\nu_{M-}(f_i^{(j)} x),a\otimes y)=S(h_i^{(j)}x,y)
+0+S(f_i^{(j)}x,ay) 
- S(x,h_i^{(j)}y)=S(f_i^{(j)}x,ay). $$
If $a$ is orthogonal to $f_i$ with respect to $S$, then (\ref{mult1}), (\ref{mult2}) 
and the result for one tensor factor give
$$S(\nu_{M-}(f_i^{(j)} x),a\otimes y)=0 +S(x,[e_i,a]^{(j)}y) +
S(f_i^{(j)}x,ay) - S(x,[e_i,a]^{(j)}y)=S(f_i^{(j)}x,ay).$$
This finishes the inductive argument. The Lemma is proved. \qed

\section{Flags, Orlik-Solomon algebra, hypergeomertic\\
 differential  forms}
  In this section we will  formulate results from \cite{SV1} which
define a map between the complex of  hypergeometric differential forms 
and the complex $C_{\bul}(\gnm^{\ast},M^{\ast})$ for a suitable Lie algebra $\gnm$ and a
$\gnm$--module $M$. 

\subsection{Complexes}\label{complexes}
 Let $W$ be an affine complex $m$--dimensional space and let $\Cc$ be a 
 configuration of hyperplanes in $W$. Define Abelian groups $\Ac^k(\Cc,\zit)$,
 $0\leq k \leq m$, as follows. $\Ac^0(\Cc,\zit)=\zit$. For $k\geq 1$, 
$\Ac^k(\Cc)$ is generated by $k$-tuples $(H_1,\ldots,H_k)$, $H_i \in \Cc$,
 subject to the relations:

 $(H_1,\ldots,H_k)=0$ if $H_1,\ldots,H_k$
 are not in general position ($\codim H_1\cap\ldots\cap H_k\ne k$);  

$(H_{\sigma(1)},\ldots,H_{\sigma(k)})=(-1)^{|\sigma|}(H_1,\ldots,H_k)$
for any permutation $\sigma \in \Sigma_k$;

$\sum_{i=1}^{k+1}(-1)^i(H_1,\ldots\Hat{H_i},\ldots,H_{k+1})=0$
for any $(k+1)$-tuple $H_1,\ldots,H_{k+1}$ which is not in general position and such that $H_1\cap\ldots\cap H_k\ne 0$.

The direct sum $\Ac^{\bullet}(\Cc,\zit)=\displaystyle\oplus_{k=0}^m
\Ac^{k}(\Cc,\zit)$ is a graded skew commutative algebra with respect to the
 multiplication $(H_1,\ldots,H_k)\cdot(H_1',\ldots,H_l')=
(H_1,\ldots,H_k,H_1',\ldots,H_l')$. \hspace{.3cm} $\Ac^{\bullet}(\Cc,\zit)$ 
 is called the {\it{Orlik-Solomon algebra}} of the configuration $\Cc$.
 
  {\bf{Flags}}. For $0\leq k\leq m$, denote by $\mathrm{Flag}^k(\Cc)$   the
 set of all flags $L^0\supset L^1\supset\cdots\supset L^k$,  where $L^i$ is an
 edge of $\Cc$ of codimension $i$. Denote by $\overline{\mathrm{Flag}}^k(\Cc)$
  the free Abelian group on $\mathrm{Flag}^k(\Cc)$ and by $Fl^k(\Cc,\zit)$ the
 quotient of  $\overline{\mathrm{Flag}}^k(\Cc)$ by the following relations. 

For every $i$, $0<i<k$, and a flag with a gap, 
$\hat{F}=(L^0\supset\cdots\supset L^{i-1}\supset L^{i+1}\supset L^k)$,
 where $L^j$ is an edge of codimension $j$, we set
 $\sum_{F\supset\hat{F}} F=0$ in  $Fl^k(\Cc,\zit)$, where the summation is
 over all flags
 $F=(\tilde{L}^0\supset\tilde{L}^k)  \in 
\mathrm{Flag}^k(\Cc)$ such that 
$\tilde{L}^j= L^j$ for all $j\ne i$.

 To define the relation between $\Ac^k(\Cc,\zit)$ and $Fl^k(\Cc,\zit)$ we 
define the following map. For $(H_1\ldots, H_k)$ in the general position, 
$H_i\in\Cc$, define $F(H_1\ldots, H_k)=(H_1\supset H_{12}\supset\cdots
\supset H_{12\ldots k})$ $\in\mathrm{Flag}^k(\Cc)$, where 
$H_{12\ldots i}=H_1\cap H_2\cap\ldots\cap H_i$. For a flag 
$F\in\mathrm{Flag}^k(\Cc)$, define a functional 
$\delta_F\in Fl^k(\Cc,\zit)^{\ast}$ as $\delta_F(F')=1$ if $F'=F$ and 
$\delta_F(F')=0$ otherwise. For $(H_1,\ldots,H_k)$ in general position, define
 a map
\begin{equation}\label{plane2flag*}
\varphi^k(H_1,\ldots,H_k)=\sum_{\sigma\in\Sigma_k}(-1)^{|\sigma|}
\delta_{F(H_{\sigma_1},\ldots,H_{\sigma_k})}.
\end{equation}
Thus  we have a homomorphism $\varphi^k:\Ac^p(\Cc,\zit)\rightarrow
Fl^p(\Cc,\zit)^{\ast}$. The following statements are from \cite{SV1}.
All groups $Fl^p(\Cc,\zit)$ are free over $\zit$.  $\Ac^p(\Cc,\zit)$ and
 $Fl^p(\Cc,\zit)$ are dual and the map $\varphi^k$ is an isomorphism.

Set $\Ac^k(\Cc)=\Ac^k(\Cc,\zit)\otimes_{\zit}\cit$ and 
$Fl^k(\Cc)=Fl^k(\Cc,\zit)\otimes_{\zit}\cit$ for all $k$.

 From now on we assume that the configuration $\Cc$ is weighted,
 that is, to any 
hyperplane $H\in\Cc$ its weight, a number $a(H)\in\cit$, is assigned. Define the
 quasiclassical weight of any edge $L$ of $\Cc$ as the sum of the weights of
 all hyperplanes that contain the edge.

 Say that a $k$-tuple $\bar{H}=(H_1,\ldots,H_k)$, $H_i\in\Cc$, is 
adjacent to a flag $F$ if there exists $\sigma \in\Sigma_k$ such that
$F=F(H_{\sigma_1},\ldots,H_{\sigma_k})$. This permutation $\sigma$ is unique.
Denote it by $\sigma(\bar{H},F)$.

Define a symmetric bilinear form $S^k$ on $Fl^k(\Cc)$. For $F,F'\in
\mathrm{Flag}^k(\Cc)$, set 
\begin{equation}\label{S-form}
S^k(F,F')=\frac{1}{k\!}\sum(-1)^{\sigma(\bar{H},F)\sigma(\bar{H},F')}
a(H_1)\ldots a(H_k),
\end{equation}
where the summation is over all $\bar{H}=(H_1,\ldots,H_k)$ adjacent to both 
$F$ and $F'$.

  The form $S^k$ is called the {\it{quasiclassical contravariant form}} of the
configuration $\Cc$. It defines a bilinear symmetric form on $FL^k(\Cc)$.
See \cite{SV1}.

 {\bf{Flag complex.}} Define a differential $d:Fl^k\rightarrow Fl^{k+1}$ 
by $d(L^0\supset\cdots\supset L^k)=\sum_{L^{k+1}} (L^0\supset\cdots\supset
 L^k\supset L^{k+1})$, where the sum is taken over all edges $L^{k+1}$ of
 codimension $k+1$ such that $L^k\supset L^{k+1}$. From the definition of 
the groups $Fl^k$ it follows that $d^2=0$.

{\bf{A complex $(\Ac^{\bullet}, d(a))$.}} Set $\omega=\omega(a)=\sum_{H\in\Cc}
a(H)H,\qquad\omega(a)\in\Ac^1$.  Define a differential 
$d=d(a):\Ac^k\rightarrow\Ac^{k+1}$ by the rule $dx=\omega(a)\cdot x$.
It is clear that $d^2=0$.

For any $k$, the quasiclassical bilinear form on $\Cc$ defines a homomorphism
\begin{equation}
S^k : Fl^k\rightarrow (Fl^k)^{\ast} \simeq \Ac^k,
\end{equation}
where $S^k(F)=(-1)^{k(k-1)/2}S(F,\, .)$.
\begin{lemma}\label{S-hom}
$S^{\bullet}$ defines a map of complexes $S^{\bullet}=S^{\bullet}(a):
(Fl^{\bullet}(\Cc),d)\rightarrow (\Ac^{\bullet}(\Cc), d(a))$.
\end{lemma}
{\bf{Note.}} There is a misprint in \cite{SV1} in the definition of $S^k$
 where the factor $(-1)^{k(k-1)/2}$ is missing.\\
{\it{Proof.}} For any edge $L$, set 
$S(L)=\sum_{H\in\Cc,\, L\subset H} a(H)H,$ \hspace{.5cm} $S(L)\in \Ac^1$.
It is easy to see that the homomorphism $S^k$ is defined by
$S^k(L^0\supset\ldots\supset L^k)=(-1)^{k(k-1)/2} S(L^1)\cdot S(L^2)\cdots
 S(L^k)$. In 
other words $S^k(L^0\supset\ldots\supset L^k)=(-1)^{k(k-1)/2}\sum a(H_1)
\ldots a(H_k) (H_1,\ldots,H_k)$, 
where the sum is over all $k$-tuples $(H_1,\ldots,H_k)$ 
such that $H_i\supset L^i$ for all $i$. Therefore, we have
\begin{eqnarray*}
\lefteqn{S^{k+1}d(L^0\supset\ldots\supset L^k)=S^{k+1}\left(
\sum_{L^{k+1},\, L^{k+1}\subset L^k}
 (L^0\supset\ldots\supset L^k\supset L^{k+1})\right)}\\
& = & (-1)^{\frac{(k+1)k}{2}}\sum a(H_1)\ldots a(H_k)a(H_{k+1})
(H_1,\ldots,H_k,H_{k+1})\\
& = &  \left( (-1)^{\frac{k(k-1)}{2}}\sum a(H_1)\ldots a(H_k)
(H_1,\ldots,H_k)\right)\cdot(-1)^k\sum a(H_{k+1})H_{k+1}\\
& = & \left( (-1)^{\frac{k(k-1)}{2}}\sum a(H_1)\ldots a(H_k)
(H_1,\ldots,H_k)\right)\cdot(-1)^k\omega(a) \\
& = & S^k(L^0\supset\ldots\supset L^k)\cdot (-1)^k\omega(a)
 = \omega(a)\cdot S^k(L^0\supset\ldots\supset L^k)
= d(a)S^k(L^0\supset\ldots\supset L^k).
\end{eqnarray*}
The second, third and fourth sum are over all $H_i$, such that $H_i\supset L^i$, for 
$1\leq i\leq k$ and $H_{k+1}\cap L^k\ne 0$. Note that  
$H_{k+1}\cap L^k = 0$ implies $(H_1,\ldots,H_k,H_{k+1})=0$, and thus the 
fourth equality is justified. The sixth one comes from the skew symmetry in
$\Ac^{\bullet}(\Cc)$.\qed
\hfill\\
Recall that we have a weighted configuration of hyperplanes in a complex $m$-
dimensional space $W$, and $a=\{a(H)|H\in\Cc\}$ are the weights. Fix an affine
 equation $l_H=0$ for each hyperplane $H\in\Cc$. Set $Y=W-\bigcup_{H\in\Cc}H$.
Consider the trivial line bundle $\Lc(a)$ over $Y$ with an integrable 
connection $d(a):\Oc\rightarrow\Omega^1$ given by
$d +\Omega (a)=d+\sum_{H\in\Cc}a(H)d\log l_H$, where $d$ is the de Rham
differential. Denote by $\Omega^{\bullet}(\Lc(a))$ the complex of $Y$-sections
of the homomorphic  de Rham complex of $\Lc(a)$.

 To any $H\in\Cc$ assign the one-form $i(H)=d\log l_H\in\Omega^1(\Lc(a))$.
This construction defines a monomorphism
$i(a):(\Ac^{\bullet}(\Cc),d(a))\rightarrow (\Omega^{\bullet}(\Lc(a)),d(a))$.
The image of this monomorphism is called the {\it{complex of the hypergeometric differential forms of weight a}}. It is denoted by $(\Ac^{\bullet}(\Cc,a),d(a))$.
The image of the homomorphism 
$i(a)S:(Fl^{\bullet}(\Cc),d)\rightarrow (\Omega^{\bullet}(\Lc(a)),d(a))$
is called the {\it{complex of the flag hypergeometric differential forms of   
 weight a}}. It is denoted by$(Fl^{\bullet}(\Cc,a),d(a))$. For further details see \cite{SV1}.

\subsection{Discriminantal configurations}
 Let $W$ be an affine complex space of dimension $m$. Let $z_1, \ldots, z_n$ be
pairwise distinct complex numbers. Denote by $\Cc_m$ a configuration in $W$
consisting of hyperplanes $H_{kl}:t_k-t_l=0;$\hspace{.5cm}$1\leq k<l\leq m$.
So $\Cc_1=\varnothing$, and $Y(\Cc_m)$ is the space of $m$-tuples of ordered
distinct 
points in $\cit$. Denote by $\Cc_{n;m}(z)$ a configuration in $W$ consisting
of hyperplanes $H_k^j: t_k-z_j=0$,\hspace{.5cm}$1\leq k\leq m$,
\hspace{.5cm}$1\leq j\leq n$, and $H_{kl},\quad 1\leq k<l\leq m$. Thus, 
$Y(\Cc_{n;m}(z))=p^{-1}(z)$ where $p:Y(\Cc_{n+m})\rightarrow Y(\Cc_n)$ is
the projection on the first $n$ coordinates. Define $\Cc_{0;m}=\Cc_m$.

{\it Edges and flags of $\Cc_{n,m}$.} For every non-empty subset 
$J=\{j_1,\ldots, j_k\}\subset [m]$ set $L_J=H_{j_1j_2}\cap H_{j_2j_3}\cap\ldots
\cap H_{j_{k-1}j_k}\in \Cc_{n;m}^{k-1}$. $L$ is an edge of codimension
 $k-1$. In particular set $L_J=W$, for $k=1$.
 For $i\in [n]$ define $L_J^i=H_{j_1}^i\cap H_{j_2}^i\cap\ldots
\cap H_{j_k}^i \in \Cc_{n;m}^k$. $L^i_J$ is an edge of codimension $k$.
 Set  $L_{\varnothing}^i=W$. 
Given non-intersecting subsets $J_1,\ldots,J_k;I_1,\ldots,I_n\subset [m]$, 
define $L_{J_1,\ldots,J_k;I_1,\ldots,I_n}=
(\bigcap_{j=1}^kL_{J_j})\cap(\bigcap_{i=1}^n L_{I_i}^i)$.

{\it Multiplication of flags.} Given two subsets $J\subset [m]$ and $I\subset [n]$,
denote by $\Cc_{J;I}\subset\Cc_{m;n}$ the subset consisting of all hyperplanes 
$H_{j_1j_2}$ with $j_1, j_2 \in J$ and $H_j^i$ with $j\in J, i\in I$. Given subsets
$J, J'\subset [m]$; $I,I'\subset[n]$ such that $J\cap J'=\varnothing$;
$I\cap I'=\varnothing$, define maps 
$\circ : \mathrm{Flag}^k(\Cc_{I;J})\times \mathrm{Flag}^l(\Cc_{I';J'})
\rightarrow \mathrm{Flag}^{k+l}(\Cc_{I\cup I',J\cup J'})$ 
as follows. For $F=F(H_1,\ldots, H_k) \in  \mathrm{Flag}^k(\Cc_{I;J})$,
 $F=F(H'_1,\ldots, H'_l) \in \mathrm{Flag}^l(\Cc_{I';J'})$, set
$F\circ F'=(H_1,\ldots, H_k,H'_1,\ldots, H'_l)$. The following lemma, 
\cite{SV1} Lemma~5.7.2, takes place .
\begin{lemma}\label{F-mult} The above map correctly defines the map
$Fl^k(\Cc_{I;J})\otimes Fl^l(\Cc_{I';J'})
\rightarrow Fl^{k+l}(\Cc_{I\cup I',J\cup J'})$. Moreover, for all 
$x\in Fl^k(\Cc_{I;J})$, $y\in Fl^l(\Cc_{I';J'})$ we have 
$x\circ y =(-1)^{kl}y\circ x$.
\end{lemma}

\subsection{Two maps of complexes}\label{two-maps}
Let $\gg$ be a Kac-Moody Lie algebra without Serre's  relations. Let 
$M=M(\Lambda_1)\otimes\cdots\otimes M(\Lambda_n)$ be a 
tensor product of Verma modules with weights 
$\Lambda_1,\ldots,\Lambda_n\in\gh^{\ast}$.
Set $\lambda=(\underbrace{1,1,\ldots,1}_{m})$. In this case the number of 
generators $(f_j)_1^r$ of $\gnm$ equals $m$, i.e. $r=m$.
Two maps of complexes
$\psi_{\bullet}$ and $\eta_{\bullet}$ are described in \cite{SV1}:
\begin{equation}\label{fi_and_psi}
\psi_p: C_p(\gnm,M)_{\lambda}\rightarrow Fl^{m-p},\qquad
\eta_p=\varphi^{-1}\circ(\psi_p^{\ast})^{-1}:C_p(\gnm^{\ast},M^{\ast})_{\lambda}
\rightarrow A^{m-p},
\end{equation}
where $\varphi$ is the map (\ref{plane2flag*}).

{\it Note.} The maps $\psi_p$ define  isomorphism of complexes. In view of Lemma~\ref{S-hom}
and \cite{SV1}, Theorem~6.6 the maps $(-1)^{(m-p)(m-p-1)/2}\eta_p$ define isomorphism of
complexes.

We will recall the explicit description of $\psi_{\bullet}$ under the above 
assumption on $\lambda$. Let $g\in\gnm$. A {\it length} $l=l(g)$ of a 
commutator $g$ is given via an inductive definition. Set $l(f_j)=1$ for 
$j=1,\ldots,m$. If $g=[g_1,g_2]$ and $l_1=l(g_1),\quad l_2=l(g_2)$, then
set $l(g)=l_1+l_2$. So $l(g)=$"the number of $f$'s in $g$''.

To every commutator $g$ assign a {\it bracket sign} $b(g)\in\zit/2\zit$ 
as follows. Set $b(f_j)=0$; $b([g_1,g_2])=b(g_1)+b(g_2)+l(g_1) \mod 2$.

To every commutator $g$ assign a flag $Fl(g)\in Fl^{l(g)-1}(\Cc_{0;|g|})$ 
as follows. Set $Fl(f_j)=\square$. If $g=[g_2,g_1]$, set $Fl(g)$ equal to 
$Fl(g_1)\circ Fl(g_2)$ completed by the edge $L_{|g|}$.

 Finally, for a commutator $g$ set $F(g)=(-1)^{b(g)}Fl(g)\in 
Fl^{l(g)-1}(\Cc_{0;|g|})$. For $I=(i_1,\ldots,i_l)\subset \{1,\ldots,m\}$ and
$1\leq i\leq n$, set $f_I=f_{i_l}\ldots f_{i_1}\in U(\gnm)$ and
$F^i(f_I)=F(H_{i_1}^i,\ldots,H_{i_l}^i)\in Fl^l(\Cc_{\{i\};I})$. 
Let $z\in C_p(\gnm,U(\gnm)^{\otimes n})_{\lambda}$ and
$z=g_p\wedge g_{p-1}\wedge\cdots g_1\otimes f_{I_n}\otimes f_{I_{n-1}}
\otimes \cdots\otimes f_{I_1}$, where all $g_i$ are commutators, $l_i=l(g_i)$.
Let $\{f_{i_1},\ldots,f_{i_m}\}$ be the list of $f_i$'s in $z$ read from 
right to left. Define $\sigma(z)\in\Sigma_m$ by $\sigma(z)(j)=i_j$. Set
\begin{equation}\label{psimap}
\psi_p(z)=(-1)^{|\sigma(z)|+\sum_{i=1}^p(i-1)(l_i-1)}
F^1(f_{I_1})\circ\ldots\circ F^n(f_{I_n})\circ F(g_1)\circ\ldots\circ F(g_p).
\end{equation}
{\bf Note.} There is a correction of the sign in the definition of
$\psi$ compared with \cite{SV1}.

{\bf Examples.}  Let $n=1$. $\psi(f_m\ldots f_1) = F(H_1^1,\ldots H_m^1)$, and
$\eta((f_{\sigma_1}\ldots f_{\sigma_m})^{\ast}) = 
(-1)^{|\sigma|}H_{\sigma_1,\sigma_2}\circ\cdots\circ H_{\sigma_{m-1},\sigma_m}
\circ H^1_{\sigma_m}$.
Compose the inclusion map 
$i(a):(\Ac^{\bullet}(\Cc),d(a))\rightarrow (\Omega^{\bullet}(\Lc(a)),d(a))$ 
(see Section~\ref{complexes}) with the map $\eta$ to get
$$i(a)\circ\eta((f_{\sigma_1}\ldots f_{\sigma_m})^{\ast})=
(-1)^{|\sigma|}d\ln(t_{\sigma_1}-t_{\sigma_2})\wedge\cdots\wedge
d\ln(t_{\sigma_{m-1}}-t_{\sigma_m})\wedge d\ln(t_{\sigma_m}-z_1).$$

 Let $I\in P(\lambda,n)$ and $ I = (i^1_{1},\ldots,i_{s_1}^1,\ldots,
i^n_{1},\ldots,i^n_{s_n})$. Since $\lambda=(1,1,\ldots,1)$, $I\in\Sigma_m$. Let
$I_j=(i^j_{s_j},\ldots,i^j_1)$ for $1\leq j\leq n$. We have
$i(a)\circ\eta ((f_{I_n})^{\ast}\otimes\cdots\otimes(f_{I_1})^{\ast})= 
(-1)^{|I|}\omega_I$. Therefore $i(a)\circ\eta(\sum_{I\in P(\lambda,n)}
(f_I)^{\ast}f_I) = \omega(z,t)$, see Section~\ref{HG-sol}.

 Let $I=(i^1_{1},\ldots,i_{s_1}^1,\ldots, i^n_{1},\ldots,i^n_{s_n})\in P(\lambda,n)$,
and $1\leq k\leq s_j$. Define
\begin{align}
   f_{I;i_k^j} &= f_{I_n}\otimes
\cdots\otimes f_{I_{j+1}}\otimes f_{i_{k+1}^j}\ldots  f_{i_{s_j}^j}\otimes
f_{I_{j-1}}\otimes\cdots\otimes f_{I_1},\notag\\
   \theta_{I;i_k^j}&= \omega_{i_1^1,\ldots,i_{s_1}^1}
\wedge\cdots\wedge \omega_{i_1^{j-1},\ldots,i_{s_{j-1}}^{j-1}}\wedge
\omega_{i_{k+1}^{j},\ldots,i_{s_j^j}}\wedge\omega_{i_1^{j+1},\ldots,i_{s_{j+1}}^{j+1}}
\wedge\cdots\wedge\omega_{i_1^n,\ldots,i_{s_n}^n}\wedge\notag\\
&\wedge [d \ln(t_{i_1^j}-t_{i_2^j})\wedge\cdots\wedge 
d\ln(t_{i_{k-1}^j}-t_{i_{k}^j})]\label{thetaform}\\
H_{I;i_k^j} &=  H^1_{i_1^1,\ldots,i_{s_1}^1}
\circ\cdots\circ H^{j-1}_{i_1^{j-1},\ldots,i_{s_{j-1}}^{j-1}}\circ
H^j_{i_{k+1}^{j},\ldots,i_{s_j^j}}\circ H^{j+1}_{i_1^{j+1},\ldots,i_{s_{j+1}}^{j+1}}
\circ\cdots\circ H^n_{i_1^n,\ldots,i_{s_n}^n}\circ\notag\\
&\circ [H_{i_1^j,i_2^j}\circ\cdots\circ H_{i_{k-1}^j,i_{k}^j}],\notag
\end{align}
where $H^p_{i_1,\ldots,i_l}=H_{i_1,i_2}\circ\cdots\circ H_{i_{l-1},i_l}\circ
H_{i_l}^p$.
\begin{lemma} Let $I\in P(\lambda,n)$. Let 
$\epsilon_k^j=k((s_j-k)+s_{j+1}+\cdots+s_n)$. Then
\begin{equation}
 i(a)\circ\eta (\Delta_{\sigma_{i_1^j},\ldots,\sigma_{i_k^j}}\otimes 
(f_{I;i_k^j})^{\ast})= (-1)^{|I|+\epsilon_k^j}\theta_{I;i_k^j}
\end{equation}  
\end{lemma}\label{g-correspondence}
{\it Proof.} The statement of the Lemma is equivalent to the equation:
\begin{align}
 \Delta_{\sigma_{i_1^j},\ldots,\sigma_{i_k^j}}\otimes 
(f_{I;i_k^j})^{\ast}&=(-1)^{|I|+\epsilon_k^j}\eta^{-1}( H_{I,i_k^j}).
\label{g-string}
\end{align}
 It is sufficient to compute the two sides on  elements of  type
$g\otimes f_{I,i_k^j}$ where $g$ is a commutator of length $k$ on 
$f_{i_1^j},\ldots,f_{i_k^j}$. Let $\sigma\in \Sigma_k$ and 
$f_{i_{\sigma_1}^j},\ldots,f_{i^j_{\sigma_k}}$ be the list of $f_{i^j}$'s 
entering $g$ from 
right to left. The left hand side  and the right hand side of (\ref{g-string}) 
evaluated on $g\otimes f_{I,i_k^j}$ give
\begin{align}
\Delta_{\sigma_{i_1^j},\ldots,\sigma_{i_k^j}}\otimes (f_{I;i_k^j})^{\ast}
(g\otimes f_{I;i_k^j})&= \Delta_{\sigma_{i_1^j},\ldots,\sigma_{i_k^j}}(g),\notag\\
(-1)^{|I|+\epsilon_k^j}\eta^{-1}( H_{I,i_k^j})(g\otimes f_{I;i_k^j})&=
(-1)^{|I|+\epsilon_k^j} \varphi(H_{I,i_k^j})(\psi(g\otimes f_{I;i_k^j})),\label{rhslema}
\end{align}
respectively. See formula (\ref{fi_and_psi}). Use the definition of $\psi$ to obtain
\begin{align}
 \psi(g\otimes f_{I;i_k^j})&= (-1)^{|\tau|}F(H^1_{i_{s_1}^1},\ldots,H^1_{i_1^1})\circ
\cdots\circ F(H^j_{i_{s_j}^1},\ldots,H^j_{i_{k+1}^j})\circ\cdots\notag\\
&\cdots\circ F(H^n_{i_{s_n}^n},\ldots,H^n_{i_1^n})\circ F(g),\notag
\end{align}
where 
$\tau=\left(\begin{smallmatrix}
  1 &. &.&.&.&.&.&.&.&.&.&.&.& m\\
 i_{s_1}^1 & \cdots & i_1^1 &\cdots &
 i_{s_j}^j & \cdots & i_{k+1}^j &\cdots & 
 i_{s_n}^n & \cdots & i_1^n &
 i_{\sigma_1}^j & \cdots & i_{\sigma_k}^j
 \end{smallmatrix}\right)$.
$I$ as an element of $\Sigma_m$ has the form
$I=\left(\begin{smallmatrix}
  1 &. &.&.&.&.&.&.&.&.& m\\
 i_{1}^1 & \cdots & i_{s_1}^1 &\cdots &
 i_{1}^j & \cdots & i_{s_j}^j &\cdots & 
 i_{1}^n & \cdots & i_{s_n}^n 
\end{smallmatrix}\right)$. Thus
\begin{equation}\label{I-tau}
|I|=(|\tau|+S_k^j+ k((s_j-k)+s_{j+1}+\cdots+s_n) + |\sigma|) \mod 2,
 \end{equation}
where $S_k^j= \sum_{l=1,\, l\ne j}^n\displaystyle\frac{s_l(s_l-1)}{2} + 
\displaystyle\frac{(s_j-k)(s_j-k-1)}{2}$.
Note that $H^p_{i_1,\ldots,i_l} =(-1)^{l(l-1)/2}H_{i_l}^p\circ H_{i_{l-1},i_l}
\circ\cdots\circ H_{i_1,i_2}$. Use the definition of $\varphi$, (\ref{plane2flag*}),
to compute
\begin{align}
\varphi(H^p_{i_1,\ldots,i_l})&=(-1)^{l(l-1)/2}
\varphi( H_{i_l}^p\circ H_{i_{l-1},i_l}\circ\cdots\circ H_{i_1,i_2})\notag\\
&=(-1)^{l(l-1)/2}\delta_{F(H_{i_l}^p,\ldots,H_{i_1}^p)} 
+ \mbox{ other } \delta-\mbox{summands}.\label{deltasum1}
\end{align}
\begin{align}
\varphi(H_{I,i_k^j})&=(-1)^{S_k^j}\varphi(
H_{i_{s_1}^1}^1\circ H_{i^1_{s_1-1},i^1_{s_1}}\circ\cdots\circ H_{i^1_1,i^1_2}
\circ\cdots\circ H_{i_{s_j}^j}^j\circ H_{i^j_{s_j-1},i^j_{s_j}}
\circ\cdots\circ H_{i^j_{k+1},i^j_{k+2}}\circ\notag\\
&\cdots\circ H_{i_{s_n}^n}^n\circ H_{i^n_{s_n-1},i^n_{s_n}}
\circ\cdots\circ H_{i^n_{1},i^n_{2}}
\circ [H_{i^j_1,i^j_2}\circ\cdots\circ H_{i^j_{k-1},i^j_k}]).\label{deltasum2}
\end{align}
 We use formulae (\ref{I-tau},\ref{deltasum1},\ref{deltasum2}) to simplify 
(\ref{rhslema}).
\begin{equation}
(-1)^{|I|+\epsilon_k^j}\eta^{-1}( H_{I,i_k^j})(g\otimes f_{I;i_k^j})=
(-1)^{|\sigma|}\varphi(H_{i^j_1,i^j_2}\circ\cdots\circ H_{i^j_{k-1},i^j_k})(F(g))
\end{equation}
The proof of Lemma~\ref{g-correspondence} is finished modulo the following result.\qed
\begin{lemma} 
Let $\eta_{\Ic}$, $\psi_{\Ic}$ be the combinatorial maps 
(\ref{fi_and_psi}), (\ref{psimap}) defined on the set of distinct indices $\Ic$,
$\Ic=\{i_1,\ldots,i_k\}\subset\{1,\ldots,m\}$. 
 Then $\eta_{\Ic} (\Delta_{i_1,\ldots,i_k})=H_{i_1,i_2}\circ\cdots\circ H_{i_{k-1},i_k}$, 
for any $k=2,\ldots,m$. 
\end{lemma}
{\it Proof.} Induction by $k$. For $k=2$, $g=[f_{i_1},f_{i_2}]$ forms a base of 
the commutators
of length 2 on $f_{i_1}$ and $f_{i_2}$. $\Delta_{i_1,i_2}(g)=1$. Since $b(g)=1$ and 
$F(g)=(-1)^{b(g)}F(H_{i_1,i_2})$ and $\sigma=\left(\begin{smallmatrix} 1 & 2 \\ 
2 & 1 \end{smallmatrix}\right)$,
we have $\eta^{-1}(H_{i_1,i_2})([f_{i_1},f_{i_2}])=\varphi(H_{i_1,i_2})(\psi(g))=
\delta_{F(H_{i_1,i_2})}((-1)^{|\sigma| + b(g)}F(H_{i_1,i_2}))=1$.

 Let $2<k\leq m$. Assume that for any $j$, $2\leq j < k$, and
$1\leq s_1<\cdots<s_j\leq k$ we have
$\eta (\Delta_{i_{s_1},\ldots,i_{s_j}})=H_{i_{s_1},i_{s_{2}}}\circ\cdots\circ 
H_{i_{s_{j-1}},i_{s_j}}$. Let $g$ be a commutator
of length $k$ on $f_{i_1},\ldots,f_{i_k}$. Then $g=[g_1,g_2]$ with $l(g_1)=l_1$, 
and $l(g_2)=l_2$, and $l_1 + l_2=k$. Let $\sigma\in\Sigma_k$ be such that
$f_{i_{\sigma_1}},\ldots,f_{i_{\sigma_k}}$ is the list of $f_i$'s in $g$ read from right 
to left.
In order to evaluate $\eta^{-1}(H_{i_1,i_2}\circ\cdots\circ H_{i_{k-1},i_k})(g)=
\varphi(H_{i_1,i_2}\circ\cdots\circ H_{i_{k-1},i_k})(\psi(g))$ remark that
\begin{align}
&\psi(g)=(-1)^{|\sigma|+b(g)}(Fl(g_2)\circ Fl(g_1),L_{|g|})=
       (-1)^{|\sigma|+l(g_1)}(F(g_2)\circ F(g_1),L_{i_1,\ldots,i_k})\notag\\
&\eta^{-1}(H_{i_1,i_2}\circ\cdots\circ H_{i_{k-1},i_k})(\psi(g))=\sum_{\tau\in\Sigma_{k-1}}
(-1)^{|\tau|}\delta_{F(H_{i_{\tau_1},i_{\tau_1+1}},\ldots, 
H_{i_{\tau_{k-1}},i_{\tau_{k-1}+1}})}(\psi(g))\label{step1}
\end{align}
Since a link corresponding to a hyperplane $H_{j,j+1}$ connects only neighbouring 
indices in a flag of a type $F(H_{i_{\tau_1},i_{\tau_1+1}},\ldots, 
H_{i_{\tau_{k-1}},i_{\tau_{k-1+1}}})$ we have
\begin{eqnarray}
\lefteqn{F(H_{i_{\tau_1},i_{\tau_1+1}},\ldots,
 H_{i_{\tau_{k-1}},i_{\tau_{k-1}+1}})=}\notag\\
&=& (\square,\cdots,(t_{i_1}=\cdots=t_{i_{\tau_{k-1}}};
 t_{i_{\tau_{k-1}+1}}=\cdots=t_{i_k}),L_{i_1,\ldots,i_k}=(t_{i_1}=\cdots=t_{i_k})).
\label{step2}\end{eqnarray}
Let $\delta_{F(H_{i_{\tau_1},i_{\tau_1+1}},\ldots, 
H_{i_{\tau_{k-}},i_{\tau_{k-1}+1}})}(\psi(g))\ne 0$. Since $(-1)^{|\sigma|+b(g)}\psi(g)=
(Fl(g_2)\circ Fl(g_1),L_{|g|})=(\square,\cdots,L_{|g_2|}\cap L_{|g_1|},L_{|g|})$
formula (\ref{step1}) implies either $L_{|g_1|}=(t_{i_1}=\cdots=t_{i_{\tau_{k-1}}})$; 
$L_{|g_2|}=(t_{i_{\tau_{k-1}+1}}=\cdots=t_{i_k})$, or 
$L_{|g_2|}=(t_{i_1}=\cdots=t_{i_{\tau_{k-1}}})$; 
$L_{|g_1|}=(t_{i_{\tau_{k-1}+1}}=\cdots=t_{i_k})$. Without loss of generality we will 
assume that 
the second case takes place, i.e. $L_{|g_2|}=(t_{i_1}=t_{i_2}=\cdots=t_{i_{\tau_{k-1}}})$; 
$L_{|g_2|}=(t_{i_{\tau_{k-1}+1}}=\cdots=t_{i_{k-1}}=t_{i_k})$. Compare the lengths of 
the flags to conclude that $\tau_{k-1}=l_2$. In order to have non-zero multiples in the 
product 
\begin{eqnarray}
  \lefteqn{\delta_{F(H_{i_{\tau_1},i_{\tau_1+1}},\ldots, 
 H_{i_{\tau_{k-1}},i_{\tau_{k-1}+1}})}(Fl(g_2)\circ Fl(g_1),L_{|g|})=}\notag\\
&=&\delta_{F(H_{i_{\tau_1},i_{\tau_1+1}},\ldots,
H_{i_{\tau_{l_2-1}},i_{\tau_{l_2-1}+1}})}(Fl(g_2))
\delta_{F(H_{i_{\tau_{l_2}},i_{\tau_{l_2}+1}},\ldots,H_{i_{\tau_{k-2}},
i_{\tau_{k-2}+1}})}(Fl(g_1))
\end{eqnarray}
 we need  $(\tau_1,\ldots,\tau_{l_2-1})$ to be a permutation of the set $(1,\ldots,l_2-1)$ and
$(\tau_{l_2},\ldots,\tau_{k-2})$ to be a permutation of the set $(l_2+1,\ldots,k-1)$.
Set 
$\tau'=\left(\begin{smallmatrix}
  1 &\cdots& l_2-1\\
  \tau_1 & \cdots & \tau_{l_2-1}
\end{smallmatrix}\right)$ and
$\tau''=\left(\begin{smallmatrix}
  1 &\cdots& l_1-1\\
  \tau_{l_2}-l_2 & \cdots & \tau_{k-2}-l_2
\end{smallmatrix}\right)$. Then $(-1)^{|\tau|}=(-1)^{|\tau'|+|\tau''|+l_1-1}$, 
$ b(g)+l_1= b(g_1)+b(g_2) \mod 2$, and
\begin{align}
&\eta^{-1}(H_{i_1,i_2}\circ\cdots\circ H_{i_{k-1},i_k})(\psi(g))=\notag\\
&=(-1)^{|\sigma|+b(g)+l_1-1}
\sum_{\tau'\in\Sigma_{l_2-1}}(-1)^{|\tau'|}
\delta_{F(H_{i_{\tau'_1},i_{\tau'_1+1}},\ldots,
H_{i_{\tau'_{l_2-1}},i_{\tau'_{l_2-1}+1}})}(Fl(g_2))\times\notag\\
&\times\sum_{\tau''\in\Sigma_{l_1-1}}(-1)^{|\tau''|}
\delta_{F(H_{i_{\tau''_1+l_2},i_{\tau''_1+l_2+1}},\ldots,
H_{i_{\tau''_{l_1-1}+l_2},i_{\tau''_{l_1-1}+l_2+1}})}(Fl(g_1))\notag\\
&=(-1)^{|\sigma|-1}\eta^{-1}(H_{i_1,i_2}\circ\cdots\circ H_{i_{l_2-1},i_{l_2}})(F(g_2))
\eta^{-1}(H_{i_{l_2+1},i_{l_2+2}}\circ\cdots\circ H_{i_{k-1},i_{k}})(F(g_1))\notag\\
&=(-1)\eta^{-1}(H_{i_1,i_2}\circ\cdots\circ H_{i_{l_2-1},i_{l_2}})(\psi(g_2))
\eta^{-1}(H_{i_{l_2+1},i_{l_2+2}}\circ\cdots\circ H_{i_{k-1},i_k})(\psi(g_1)).\label{step3}
\end{align}
The last equality holds because
$\sigma=\left(\begin{smallmatrix}
  1 &\cdots& k\\
  \sigma_1 & \cdots & \sigma_k
\end{smallmatrix}\right)=
\left(\begin{smallmatrix}
  1 &\cdots& l_2\\
  \sigma_1 & \cdots & \sigma_{l_2}
\end{smallmatrix}\right)\times
\left(\begin{smallmatrix}
  l_2+1 &\cdots& l_k\\
  \sigma_{l_2+1} & \cdots & \sigma_{l_k}
\end{smallmatrix}\right).$ Using the inductive hypothesis rewrite (\ref{step3}) as
\begin{align}
\eta^{-1}(H_{i_1,i_2}\circ\cdots\circ H_{i_{k-1},i_k})(\psi(g))&=
(-1)\Delta_{i_1,\ldots,i_{l_2}}(g_2)\Delta_{i_{l_2+1},\ldots,i_k}(g_1)\notag\\
&=\Delta_{i_1,\ldots,i_k}([g_1,g_2]).
\end{align}\qed

\section{Derivation of the dynamical differential equation}\label{DynDE}
  In this section $\gg$ will be a Kac-Mody Lie algebra without Serre's relations,
$\lambda=(1,1,\ldots,1)$, $r=m$. We will work in a weight space 
$M_{\lambda}$ of the module $M=M(\Lambda_1)\otimes\cdots\otimes M(\Lambda_n)$.
We will differentiate the hypergeometric form $\omega(z,t)$, express
the result in terms of the complex $C_{\bullet}(\gnm^{\ast},M^{\ast})$, and derive
the Dynamical differential equation in the form (\ref{DynKZ*}).

  The integrand of a hypergeometric solution have the following form, see 
Section~\ref{HG-sol}.
\begin{align}
\Phi_{\mu}^{\frac{1}{\kappa}}\omega&=\exp(\frac{1}{\kappa} 
(-\sum_{i=1}^m\la\alpha_{c(i)},\mu\ra t_i +
\sum_{j=1}^n\la\Lambda_j,\mu\ra z_j))\Phi^{\frac{1}{\kappa}}, \mbox{ where } \notag\\ 
\Phi(z,t)&=\prod_{i<j}(z_i-z_j)^{(\Lambda_i,\Lambda_j)}
\prod_{k,j}(t_k-z_j)^{-(\alpha_{c(k)},\Lambda_j)}
\prod_{k<l}(t_k-t_l)^{(\alpha_{c(k)},\alpha_{c(l)})}.\notag
\end{align}
Fix $\mu'\in \gh$ and let $\partial_{\mu'}$ be the partial derivative with respect to the 
parameter $\mu$ in the direction of $\mu'$. Then
\begin{equation}\label{diff1}
\kappa\partial_{\mu'}(\Phi_{\mu}^{\frac{1}{\kappa}}\omega) = 
(-\sum_{i=1}^m\la\alpha_{c(i)},\mu'\ra t_i +\sum_{j=1}^n\la\Lambda_j,\mu'\ra z_j)
\Phi_{\mu}^{\frac{1}{\kappa}}\sum_{I\in P(\lambda,n)}(-1)^{|I|}\omega_I f_Iv.
\end{equation}
 Let $I=(i_1^1,\ldots,i^1_{s_1};\ldots;i_1^n,\ldots,i^n_{s_n}) \in P(\lambda,n)$. Moreover
$I\in\Sigma_m$ because of the form of $\lambda$. Since $r=m$ we have $c(i)=i$. 
Set $t_{i_{s_j+1}^j}=z_j$ and $\alpha_{i_{k;j}}=\alpha_{i_k^j}+ \alpha_{i_{k-1}^j}+\cdots
+\alpha_{i_1^j}$. Rearrange the following expression:
\begin{eqnarray}
\lefteqn{(-\sum_{i=1}^m\la\alpha_{i},\mu'\ra t_i +\sum_{j=1}^n\la\Lambda_j,\mu'\ra z_j)
\omega_I = 
(-\sum_{j=1}^n\sum_{k=1}^{s_j}\la\alpha_{i_k^j},\mu'\ra t_{i_k^j} +
\sum_{j=1}^n\la\Lambda_j,\mu'\ra z_j)\omega_I}\notag\\
&=& (-\sum_{j=1}^n\sum_{k=1}^{s_j}\la\alpha_{i_k^j},\mu'\ra(t_{i_k^j}-t_{i_{k+1}^j} +
t_{i_{k+1}^j}-t_{i_{k+2}^j}+\cdots + t_{i_{s_j}^j} - z_j +z_j)+
\sum_{j=1}^n\la\Lambda_j,\mu'\ra z_j)\omega_I \notag\\
&=& (-\sum_{j=1}^n\sum_{k=1}^{s_j}\la\alpha_{i_{k;j}},\mu'\ra 
(t_{i_k^j}-t_{i_{k+1}^j}))\omega_I +
(\sum_{j=1}^n z_j\la\Lambda_j-\alpha_{i_1^j}-\cdots-\alpha_{i_{s_j}^j},\mu'\ra)\omega_I.
\label{diff2}\end{eqnarray}
\begin{lemma}
Let $u(\mu,z)=\int_{\gamma(z)}\Phi_{\mu}^{\frac{1}{\kappa}}\omega$. Then
\begin{equation}
\kappa\partial_{\mu'}u-\sum_{j=1}^n z_j {\mu'}^{(j)}u = \sum_{I\in P(\lambda,n)}
(-1)^{|I|}\int_{\gamma(z)}\Phi_{\mu}^{\frac{1}{\kappa}}
(-\sum_{j=1}^n\sum_{k=1}^{s_j}\la\alpha_{i_{k;j}},\mu'\ra(t_{i_k^j}-t_{i_{k+1}^j}))\omega_If_Iv
\end{equation}
\end{lemma}
{\it Proof.} Combine formulae (\ref{diff1}), (\ref{diff2}) with the fact
${\mu'}^{(j)}f_Iv=\la\Lambda_j-\alpha_{i_1^j}-\cdots-\alpha_{i_{s_j}^j},\mu'\ra f_Iv$ to obtain 
the result. \qed
\begin{lemma}\label{L-operator} Define an operator $L$ by 
$Lu=\kappa\partial_{\mu'}u-\sum_{j=1}^n z_j {\mu'}^{(j)}u$. Then
$$ Lu=
\sum_{\alpha>0}\frac{\la\alpha,\mu'\ra}{\la\alpha,\mu\ra}\Delta_{+,\alpha}u.$$
\end{lemma}
{\it Proof.} Let $I\in P(\lambda,n)$ and  $1\leq k\leq s_j$. 
\begin{align}
&(t_{i_k^j}-t_{i_{k-1}^j})\omega_I = \omega_{i_1^1,\ldots,i_{s_1}^1}\wedge\cdots\wedge
(t_{i_k^j}-t_{i_{k-1}^j})\omega_{i_1^j,\ldots,i_{s_j}^j}\wedge\cdots\wedge
\omega_{i_1^n,\ldots,i_{s_n}^n}\notag\\
&=\omega_{i_1^1,\ldots,i_{s_1}^1}\wedge\cdots\wedge
[ d\ln(t_{i_1^j}-t_{i_2^j})\wedge\cdots\wedge d\ln(t_{i_{k-1}^j}-t_{i_{k}^j})\wedge
d(t_{i_{k}^j}-t_{i_{k+1}^j})\wedge \notag\\
&\wedge d\ln(t_{i_{k+1}^j}-t_{i_{k+2}^j})\wedge\cdots\wedge d\ln(t_{i_{s_j}^j}-z_j)] 
\wedge\cdots\wedge\omega_{i_1^n,\ldots,i_{s_n}^n}\notag\\
&=(-1)^{\epsilon_k^j}\omega_{i_1^1,\ldots,i_{s_1}^1}
\wedge\cdots\wedge \omega_{i_{k+1}^j,\ldots,i_{s_j}^j}
\wedge\cdots\wedge\omega_{i_1^n,\ldots,i_{s_n}^n}\wedge\notag\\
&\wedge[d \ln(t_{i_1^j}-t_{i_2^j})\wedge\cdots\wedge 
d\ln(t_{i_{k-1}^j}-t_{i_{k}^j})]\wedge d(t_{i_{k}^j}-t_{i_{k+1}^j})\notag\\
&=(-1)^{\epsilon_k^j}\theta_{I;i_k^j}\wedge d(t_{i_{k}^j}-t_{i_{k+1}^j})=
(-1)^{\epsilon_k^j+(m-1)} d(t_{i_{k}^j}-t_{i_{k+1}^j})\wedge\theta_{I;i_k^j},
\label{diff3}
\end{align}
where $\epsilon_k^j=k(s_{j+1}+\cdots s_n+s_j-k)$. It is clear that 
$d\theta_{I;i_k^j}=0$. Thus
\begin{equation}
\kappa d_t (\Phi_{\mu}^{\frac{1}{\kappa}}\theta_{I;i_k^j}) =(-\sum_{i=1}^m
\la\alpha_i,\mu\ra dt_i+
 \Phi^{-1}(d_t\Phi))\Phi_{\mu}^{\frac{1}{\kappa}}\theta_{I;i_k^j}.
\end{equation}
 Rearrange as in formula (\ref{diff2}) and simplify to get
\begin{align}
 (\sum_{j=1}^n\sum_{k=1}^{s_j}\la\alpha_{i_{l;j}},\mu'\ra d(t_{i_l^j}-t_{i_{l+1}^j}))\wedge
\theta_{I;i_k^j}& =\Phi^{-1}(d_t\Phi))\wedge\theta_{I;i_k^j} - 
\kappa\Phi_{\mu}^{-\frac{1}{\kappa}}d_t (\Phi_{\mu}^{\frac{1}{\kappa}}\theta_{I;i_k^j}) 
\notag\\
\la\alpha_{i_{k;j}},\mu'\ra d(t_{i_k^j}-t_{i_{k+1}^j}))\wedge
\theta_{I;i_k^j}& =\Phi^{-1}(d_t\Phi))\wedge\theta_{I;i_k^j} - 
\kappa\Phi_{\mu}^{-\frac{1}{\kappa}}d_t (\Phi_{\mu}^{\frac{1}{\kappa}}\theta_{I;i_k^j})
\label{diff4}
\end{align}
 Since $\gamma(z)$ is a cycle, formulae (\ref{diff3}), (\ref{diff4}) 
allow us to rewrite $Lu$ as
\begin{align}
Lu &= \sum_{I\in P(\lambda,n)} (-1)^{|I|}\int_{\gamma(z)}(-1)^{\epsilon_k^j}
(-\sum_{j=1}^n\sum_{k=1}^{s_j}\frac{\la\alpha_{i_{k;j}},\mu'\ra}{\la\alpha_{i_{k;j}},\mu\ra}
(\Phi_{\mu}^{\frac{1}{\kappa}}\Phi^{-1}(d_t\Phi))\wedge\theta_{I;i_k^j} - 
\kappa d_t (\Phi_{\mu}^{\frac{1}{\kappa}}\theta_{I;i_k^j})\notag\\
&=\sum_{I\in P(\lambda,n)}\int_{\gamma(z)}\Phi_{\mu}^{\frac{1}{\kappa}}
(-\sum_{j=1}^n\sum_{k=1}^{s_j}\frac{\la\alpha_{i_{k;j}},\mu'\ra}{\la\alpha_{i_{k;j}},\mu\ra}
(-1)^{m-1}\frac{d_t\Phi}{\Phi}\wedge(i(a)\circ\eta(\Delta_{i_1^j,\ldots,i_k^j}
\otimes(f_{I,i_k^j})^{\ast}))
\end{align}
Since $\theta_{I,i_k^j}$ is closed, its differential in the complex of the hypergeometric 
differential forms reduces to multiplication by $i(a)(\Omega(a))=d_t\Phi/\Phi$. Taking
into account Lemma~\ref{S-hom} we have
\begin{align}
Lu&=\sum_{I\in P(\lambda,n)}f_I\int_{\gamma(z)}\Phi_{\mu}^{\frac{1}{\kappa}}
(-\sum_{j=1}^n\sum_{k=1}^{s_j}\frac{\la\alpha_{i_{k;j}},\mu'\ra}{\la\alpha_{i_{k;j}},\mu\ra}
(-1)^{m-1}d(a)(i(a)\circ\eta(\Delta_{i_1^j,\ldots,i_k^j}
\otimes(f_{I,i_k^j})^{\ast}))\notag\\
&=\sum_{I\in P(\lambda,n)}f_I\int_{\gamma(z)}\Phi_{\mu}^{\frac{1}{\kappa}}
(-\sum_{j=1}^n\sum_{k=1}^{s_j}\frac{\la\alpha_{i_{k;j}},\mu'\ra}{\la\alpha_{i_{k;j}},\mu\ra}
(i(a)\circ\eta(d\Delta_{i_1^j,\ldots,i_k^j}
\otimes(f_{I,i_k^j})^{\ast}))\notag
\end{align}
\begin{align}
&=\sum_{I\in P(\lambda,n)}f_I\int_{\gamma(z)}\Phi_{\mu}^{\frac{1}{\kappa}}
(-\sum_{j=1}^n\sum_{k=1}^{s_j}\frac{\la\alpha_{i_{k;j}},\mu'\ra}{\la\alpha_{i_{k;j}},\mu\ra}
(i(a)\circ\eta(\Delta_{i_1^j,\ldots,i_k^j}(f_{I,i_k^j})^{\ast}))\notag\\
&=\sum_{I\in P(\lambda,n)}f_I\int_{\gamma(z)}\Phi_{\mu}^{\frac{1}{\kappa}}
(-\sum_{j=1}^n\sum_{k=1}^{s_j}\frac{\la\alpha_{i_{k;j}},\mu'\ra}{\la\alpha_{i_{k;j}},\mu\ra}
(i(a)\circ\eta(\Delta_{i_1^j,\ldots,i_k^j}e_{i_k^j}\ldots e_{i_1^j}(f_I)^{\ast}))
\end{align}
 Note that $\alpha_{i_{k;j}}$, $k=1,\ldots,s_j$ describe all $\lambda$-admissible roots 
such that $\Delta_{+,\alpha}^{(j)}(f_I)^{\ast}\ne 0$, see
Section~\ref{symsection}. Therefore
\begin{align}
Lu&=-\sum_{I\in P(\lambda,n)}f_I\int_{\gamma(z)}\Phi_{\mu}^{\frac{1}{\kappa}}
(i(a)\circ\eta(\sum_{\alpha>0}\frac{\la\alpha,\mu'\ra}{\la\alpha,\mu\ra}
\sum_{(i_1,\ldots,i_{m'})\in P(\lambda_{\alpha},1)}\Delta_{i_1,\ldots,i_{m'}}e_{i_{m'}}
\ldots e_{i_1}(f_I)^{\ast}))\notag\\
&=-\sum_{I\in P(\lambda,n)}f_I\int_{\gamma(z)}\Phi_{\mu}^{\frac{1}{\kappa}}
(i(a)\circ\eta(\sum_{\alpha>0}\frac{\la\alpha,\mu'\ra}{\la\alpha,\mu\ra}
-\Delta_{-,\alpha}(f_I)^{\ast}))\notag\\
&=-\sum_{K,I\in P(\lambda,n)}f_I\int_{\gamma(z)}\Phi_{\mu}^{\frac{1}{\kappa}}
(i(a)\circ\eta(\sum_{\alpha>0}\frac{\la\alpha,\mu'\ra}{\la\alpha,\mu\ra}
\la-\Delta_{-,\alpha}(f_I)^{\ast}),f_K\ra(f_K)^{\ast})\notag\\
&=-\sum_{\alpha>0}\frac{\la\alpha,\mu'\ra}{\la\alpha,\mu\ra}
\sum_{K,I\in P(\lambda,n)}\la-\Delta_{-,\alpha}(f_I)^{\ast}),f_K\ra u_K f_I
\end{align}
\begin{align}
Lu&=-\sum_{\alpha>0}\frac{\la\alpha,\mu'\ra}{\la\alpha,\mu\ra}
\la-\Delta_{-,\alpha}(f_I)^{\ast}),u\ra f_I
=-\sum_{\alpha>0}\frac{\la\alpha,\mu'\ra}{\la\alpha,\mu\ra}\Delta_{+,\alpha}u\qquad\square
\end{align}
The statement of Lemma~\ref{L-operator} is equivalent to the Dynamical differential
equation in the direction of $\mu'$ for a function $u$ with values in
the $(1,1,\ldots,1)$ weight space of a $\gg$ module $M$. The Symmetrization Lemma~\ref{sym-lemma}
deduces the general case from this one.

\section{Main theorems}\label{main-theorems}
In this section we conclude the proofs of the the Theorems from
Section~\ref{HG-sol} 
in the setting of Kac-Moody Lie algebras without Serre's relations. Then we deduce 
the corresponding results for any simple Lie algebra.

Let $\gg$ be a Kac-Moody Lie algebra without Serre's relations. Let $\lambda\in\nit^r$. 
Let $M=M(\Lambda_1)\otimes\cdots\otimes M(\Lambda_n)$ be a tensor product of Verma
modules for $\gg$ with highest weights $\Lambda_1,\ldots,\Lambda_n\in\gh^{\ast}$. Let 
$u(\mu,z)=\sum_{I\in P(\lambda,n)}u_I f_I$ be a hypergeometric integral with values in the 
weight space $M_{\lambda}$ as described in
Section~\ref{HG-sol}, i.e. $u_I=\int_{\gamma(z)}\Phi_{\mu}^{\frac{1}{\kappa}}
(\sum_{\sigma\in \S(I)}(-1)^{|\sigma|}\omega_{I,\sigma})$.
\begin{theorem} \label{gmain1}
The function $u(\mu,z)$ solves the {\rm KZ} equations (\ref{genKZ}) in $M_{\lambda}$.
\end{theorem} 
{\it Proof.} The proof given in Section~\ref{HG-sol} holds, because all relations we used
are proved in \cite{SV1} in the general setting described above.\qed
\begin{theorem}\label{gmain2}
 The function $u(\mu,z)$ solves the dynamical differential equations
(\ref{DynKZ}) in $M_{\lambda}$.
\end{theorem}  
{\it Proof.} Lemma~\ref{sym-lemma} reduces the case of a general weight space 
$M_{\lambda}$ to the case
of a weight space $M_{\tlam}$, where $\tlam=(\underbrace{1,1,\ldots,1}_{m})$. 
Lemma~\ref{L-operator} derives the Theorem in that case.\qed

{\bf Proof of Theorem~\ref{main2}.} Combine Corollary~\ref{corCD} and Theorem~\ref{gmain2}
to derive the dynamical differential equations for any Kac-Moody Lie
algebra. In  
particular we have it for a simple Lie algebra. \qed

 Finally we will prove a determinant formula which establishes a basis of solutions for the
system of KZ and dynamical differential equations in a weight space
$M_{\lambda}$. From that
formula we will derive the compatibility of the system of KZ and Dynamical differential
equations.

Fix $\lambda\in\nit^r$. Fix a basis $(f_Iv)_{I\in P(\lambda,n)}$ of
the weight space $M_{\lambda}$. Assume that a set $(\gamma_I(z))_{I\in P(\lambda,n)}$
of horizontal families of twisted cycles in $\{z\}\times\cit^m$ is given

Denote $u_{IJ}=\int_{\gamma_I(z)}\Phi_{\mu}^{\frac{1}{\kappa}}
(\sum_{\sigma\in S(J)}(-1)^{|\sigma|}\omega_{J,\sigma})$.

\begin{proposition} \label{det-gen}
Let $\delta_{\alpha}= \tr_{M_{\lambda}}(\Delta_{+,\alpha})$ for a
positive root $\alpha$ of 
$\gg$. Denote $\epsilon_{ij}= \tr_{M_{\lambda}}(\Omega_{ij,+})$. Then we have

(a) For any  horizontal families of twisted cycles $(\gamma_I(z))_{I\in P(\lambda,n)}$
 in $\{z\}\times\cit^m$, there exists
a constant  $C=C(\Lambda_1,\ldots,\Lambda_n,\lambda,\kappa)$ such that
\begin{equation}
 \det (u_{IJ})= C\exp(\sum_{i=1}^n \frac{z_i}{\kappa}\tr_{M_{\lambda}}(\mu^{(i)}))
\prod_{\alpha>0}\la\alpha,\mu\ra^{\frac{\delta_{\alpha}}{\kappa}}
\prod_{i<j}(z_i-z_j)^{\frac{\epsilon_{ij}}{\kappa}}.
\end{equation}
In the first product only finite number of factors are different from 1, i.e.
$\delta_{\alpha}\ne 0$ if and only if $0<\alpha\leq\lambda$.

(b) For generic values of the parameters 
$(\Lambda_j)_{j=1}^m,(\alpha_i)_{i=1}^r,\kappa$ in a neighbourhood of a generic point
$(\mu,z)\in\gh\times\cit^n$ we can choose cycles $(\gamma_I(z))_{I\in P(\lambda,n)}$
such that the constant $C$ from (a) is non-zero.
 Moreover the set of functions 
$\{u^I=\sum_{J\in P_(\lambda,n)}u_{IJ}f_J\}_{I\in P(\lambda,n)}$ form a fundamental system of 
solutions for the system of KZ and dynamical differential equations. 
\end{proposition}
{\it Proof.} Part (a)  is a corollary of Theorems~\ref{gmain1} and
\ref{gmain2}. We will prove part (b) for values of the parameters such that 
all numbers $(\alpha_i,\alpha_j)/\kappa$, $-(\alpha_i,\Lambda_k)/\kappa$ 
have positive real parts for
$1\leq i,j\leq r$, $1\leq k\leq n$, and for a point $(\mu,z)$ such $z\in\rit^n$,
$z_1<z_2<\cdots<z_n$ and $\la\alpha_i,\mu\ra/\kappa>0$ for any 
$i=1,\ldots,r$.  For generic values of  
$(\Lambda_j)_{j=1}^m,(\alpha_i)_{i=1}^r, z, \mu, \kappa$
\hspace{12pt} (b) holds by analytic continuation.

 The case $\lambda=(\underbrace{1,1,\ldots,1}_m)$. Set 
$f_0(t)=\sum_{j=1}^m\la\alpha_{c_{\lambda}(j)},\mu\ra t_j$.
Let $I=(i_1^1,\ldots,i^1_{s_1};\ldots;$ $i^n_1,\ldots,i^n_{s_n})\in P(\lambda,n)$.
Set $\gamma_I(z)=\{t\in\rit^{n}\mbox{ : } z_j<t_{i^j_1}<\cdots<t_{i^j_{s_j}}<z_j+1 
\mbox{ for all } j=1,\ldots,n\}$, where $z_{n+1}=\infty$. 
Note that $\{\gamma_I(z)\}_{I\in P(\lambda,n)}$ is the set of all domains for the 
configuration of hyperplanes $H_{ij}:t_i-t_j=0$, $H_i^k: t_i-z_k=0$, $1\leq i<j\leq m$, 
$1\leq k\leq n$  which are either bounded, or the limit of $f_0$ on them is $+\infty$ when
$\| t\|\rightarrow\infty$. In \cite{Z} a linearly independent set of hypergeometric
differential $n$-forms, called $\beta\mathbf{nbc}$ differential $n$-forms,
associated to those domains is defined. An explicit non-vanishing 
formula for the corresponding determinant is given in \cite{MTV} Theorem~6.2, see also \cite{DT}. 
Since $\lambda=(1,1,\ldots,1)$ the space of hypergeometric $n$-forms is isomorphic
to the space $C_0(\gnm^{\ast},M^{\ast})_{\lambda}$, see Section~\ref{two-maps}. The latter has basis
$(f_Iv^{\ast})_{I\in P(\lambda,n)}$ which gives the basis 
$(\omega_I=i(a)\circ\eta(f_Iv^{\ast}))_{I\in P(\lambda,n)}$ of 
the space of hypergeometric $n$-forms. Since this basis and the
$\beta \mathbf{nbc}$ set have the same cardinality the non-zero determinant formula
for the integrals of $\beta \mathbf{nbc}$ forms over the domains 
$\{\gamma_I(z)\}_{I\in P(\lambda,n)}$
implies a non-zero determinant formula for the integrals of 
$(\omega_I)_{I\in P(\lambda,n)}$ over
the same domains. Since the determinant is non-zero at one point $(\mu,z)$, it is non-zero at 
any point $(\mu,z)$ under the above conditions on the parameters.

 The case of generic $\lambda\in\nit^r$. Consider 
$C_0(\gnm^{\ast},M^{\ast})_{\lambda}$ and $C_0(\tnm,\widetilde{M}^{\ast})_{(1,1,\ldots,1)}$
as in Section~\ref{symsection}. A basis for the $\Sigma_{\lambda}$-symmetric hypergeometric 
differential forms is given by $(\omega_I)_{I\in P(\lambda,n)}$, where
$\omega_I=\sum_{J\in \S(I)}\tom_J$ and $\tom_J= i(a)\circ\eta(\tf_J^{\ast})$. 
$P((\underbrace{1,1,\ldots,1}_m),n)=\displaystyle{\cup_{I\in P(\lambda,n)}}S(I)$ a disjoint
union. Thus the set $(\omega_I)_{I\in P(\lambda,n)}$ consists of  linearly independent 
forms in the space of all 
hypergeometric forms. The integral pairing described in the previous case is non-degenerate.
Therefore there there exists a subset  of the set $(\gamma_J(z))_{J\in P((1,1,\ldots,1),n)}$ 
indexed by the set $P(\lambda,n)$ such that the corresponding determinant is non-zero.
\qed
\begin{corollary} The system consisting of the union of KZ and Dynamic differential 
equations for any Kac-Moody Lie algebra with (or without)
 Serre's relations is a compatible system of differential equations.
\end{corollary}
{\bf Remark.} An algebraic proof of the compatibility of the system of KZ equations is given 
in \cite{SV1}.\\
{\it Proof.} Let us write the differential operators which determine the 
KZ equations (\ref{genKZ}) and the dynamical equations (\ref{DynKZ}) in the form
\begin{equation}
\mbox{ KZ: }\quad \frac{\partial}{\partial z_j}+ B_j,\qquad 
\mbox{ Dynamical: } \frac{\partial}{\partial \mu'} + C_{\mu'}, \quad 
\mbox{ where } j=1,\ldots,n \mbox{ and } \mu'\in\gh.
\end{equation}
The operators $B_j$ and $C_{\mu'}$ are linear for any $j=1,\ldots,n,\,\,\mu'\in\gh$.
In order to prove the compatibility of the system of KZ and Dynamical 
differential equations
we need to check $[\frac{\partial}{\partial z_j}+ B_j,
\frac{\partial}{\partial z_k} + B_k]=0$,
$[\frac{\partial}{\partial z_j}+ B_j,\frac{\partial}{\partial \mu'} + C_{\mu'}]=0$, and
$[\frac{\partial}{\partial \mu'}+ C_{\mu'},\frac{\partial}{\partial \mu''} + C_{\mu''}]=0$.

First consider the case of a Kac-Moody Lie algebra without Serre's relations,
$\gg$, acting on  a tensor product of highest weight modules $M$.
We have
\begin{equation}\label{comm-op}
[\frac{\partial}{\partial z_j}+ B_j,\frac{\partial}{\partial \mu'} + C_{\mu'}]
=(\frac{\partial}{\partial z_j}C_{\mu'})-(\frac{\partial}{\partial \mu'}B_j)
+[B_j,C_{\mu'}]
\end{equation}
The result is a linear  operator with meromorphic coefficients depending on 
parameters $\{z,\mu$, $(\alpha_i)_{i=1}^r$, $(\Lambda_j)_{j=1}^n,\kappa \}$.
Analogously, the commutators 
$[\frac{\partial}{\partial z_j}+ B_j,\frac{\partial}{\partial z_k} + B_k]$ and
$[\frac{\partial}{\partial \mu'}+ C_{\mu'},\frac{\partial}{\partial \mu''} + C_{\mu''}]$ are
linear  operators with meromorphic coefficients depending on the above set of 
parameters, where $j,k=1,\ldots,n$, $\mu',\mu''\in\gh$. 

 It is enough to show the commutativity of the above operators for generic values of 
the parameters. Then the commutators will be zero for any values of the parameters by
analytic continuation.

Take such parameters $\{z,\mu$, $(\alpha_i)_{i=1}^r$, $(\Lambda_j)_{j=1}^n,\kappa \}$
that the set of hypergeometric solutions of the system of KZ and dynamical differential 
equations forms a basis of $M$. According to
Proposition~\ref{det-gen} (b) this is  a generic 
choice of parameters. Since the KZ and the Dynamical differential
operators act as zero
on the set of hypergeometric solutions, their commutators also act as zero on the same 
set. Therefore the commutators act as zero on the $\gg$--module $M$.

  Finally, consider a Kac-Moody Lie algebra with Serre's relations
$\bar{\gg}=\gg/\ker(S:\gg\rightarrow\gg^{\ast})$ which acts on
$L=M/\ker(S:M\rightarrow M^{\ast})$. Corollary~\ref{corCD} and \cite{SV1}  Corollary~7.2.11 show
that the Dynamical and the KZ operators for $\gg$ correspond to the
 the the Dynamical and the KZ operators for $\bar{\gg}$ under this factorization. Then the 
commutativity of the operators on $\gg$ implies that they are commutative on $\bar{\gg}$ as well.
\qed

\end{document}